\documentclass[notitlepage]{scrartcl}
\usepackage[margin=1in]{geometry}
\usepackage{amsmath}
\usepackage{xcolor}
\usepackage{amssymb}
\usepackage{amsthm}
\usepackage{abstract}
\usepackage{hyperref}
\newcommand{\footremember}[2]{%
    \footnote{#2}
    \newcounter{#1}
    \setcounter{#1}{\value{footnote}}%
}
 
\title{Site Percolation on Pseudo-Random Graphs}
\author{%
Sahar Diskin \footremember{alley}{School of Mathematical Sciences, Tel Aviv University, Tel Aviv 6997801, Israel. Email: sahardiskin@mail.tau.ac.il.}%
\and Michael Krivelevich \footremember{trailer}{School of Mathematical Sciences, Tel Aviv University, Tel Aviv 6997801, Israel. Email:
krivelev@tauex.tau.ac.il. Research supported in part by USA–Israel BSF grant 2018267 and by ISF grant
1261/17.}%
}
\begin{document}
	\maketitle
\begin{onecolabstract}
    We consider vertex percolation on pseudo-random $d-$regular graphs. The previous study by the second author established the existence of phase transition from small components to a linear (in $\frac{n}{d}$) sized component, at $p=\frac{1}{d}$. In the supercritical regime, our main result recovers the sharp asymptotic of the size of the largest component, and shows that all other components are typically much smaller. Furthermore, we consider other typical properties of the largest component such as the number of edges, existence of a long cycle and expansion. In the subcritical regime, we strengthen the upper bound on the likely component size.   
\end{onecolabstract}

\section{Introduction and Main Results}
	In 1957, Broadbent and Hammersley initiated the study of percolation theory~\cite{broadbent and hammersley} in order to model the flow of fluid through a medium with randomly blocked channels. Since then, the theory of percolation has been studied extensively (see~\cite{kesten},~\cite{grimmet} and~\cite{percolation book} for systematic coverage). 
	
	Two models have been of significant interest: \textit{bond percolation} and \textit{site percolation}. In the site percolation model, given a graph $G=(V,E)$, we form a random subset of the vertices of $G$, $R\subseteq V$, and consider the subgraph induced by this subset, $G[R]$. In the bond percolation model, we form a random subgraph by choosing a random subset of the edges of $G$. While the bond percolation model has perhaps been more studied between the two, both are quite important and both serve as models of physical processes.
	
	In both models, early research focused on percolation on specific graphs, such as the lattice $\mathbb{Z}^d$ in the infinite case, and similar structures such as the $d$-dimensional hypercube and torus in the finite case. In their pioneering paper from 2004~\cite{ABS}, Alon, Benjamini and Stacey studied bond percolation on a large family of graphs --- expanders. In their paper, they established that typically there is at most one linear-sized component; they also showed the typical existence of a linear-sized component in the case of constant-degree high-girth expanders in the supercritical regime. Subsequent work of Frieze, Krivelevich and Martin (in the case of growing-degree expanders)~\cite{random subgraphs} and of Krivelevich, Lubetzky and Sudakov (in the case of constant-degree high-girth expanders)~\cite{high-girth} recovered the asymptotics and other typical properties of the giant in the supercritical regime, and discussed the typical size of the second largest component. Here we aim to achieve analogous results in the case of site percolation on an important family of expanders --- the $(n,d,\lambda)$ pseudo-random graphs.
	
	An $(n,d,\lambda)$\textit{-graph} $G$ is a $d$-regular graph on $n$ vertices, where its eigenvalues (i.e. the eigenvalues of the adjacency matrix of $G$) $d=\lambda_1\ge \lambda_2\ge\cdots\ge \lambda_n$ satisfy $\lambda=max(|\lambda_2|,|\lambda_n|)$. The greater the ratio between $d$ and $\lambda$, the more tightly the distribution of the edges of the graph approaches that of a random graph $G(n,d/n)$, which can be seen by the expander mixing lemma (stated as Lemma 2.2 in this paper). Thus, $(n,d,\lambda)$-graphs serve frequently as a standard model of pseudo-random graphs. We refer to~\cite{survey pseudo} for a comprehensive survey on the subject of pseudo-random graphs. We note that for any meaningful results, we need to require some control over the edge distribution of the base graph $G$, which we obtain by focusing on this type of graphs. For example, if $G$ is a disjoint union  of cliques of size $d+1$, all the connected components of $G[R]$ are of size at most $d+1$ deterministically.
	
	For several concrete $d$-regular base graphs, it is known that around probability $p=\frac{1}{d}$ abrupt changes occur in site percolation, in particular the emergence of a giant component whose size is proportional to $n/d$ --- for example, in the $d$-dimensional cube $\mathbb{Q}^d$~(\cite{cube 1},~\cite{cube 2}) and the $d$-dimensional Hamming torus \cite{hamming torus}. In 2016, Krivelevich established the phase transition in site percolation on $(n,d,\lambda)$-graphs~\cite{site-percolation}, showing a jump from typically logarithmic-sized components for $p=\frac{1-\epsilon}{d}$ (\textit{subcritical phase}) to a typically linear (in $n/d$) sized component for $p=\frac{1+\epsilon}{d}$ (\textit{supercritical phase)}. Since $R$ is typically of size proportional to $n/d$ given that $p$ is around $\frac{1}{d}$, it is natural to scale the obtained structures relative to this size, and to require in particular that $d=o(n)$. We aim to improve upon these results, and in particular to recover the asymptotic order of a largest component, and to show that typically there is a unique giant component, with the other components being much smaller.
	
	Trying to recover the asymptotic of a giant and to establish its uniqueness in our setting, a natural first approach is to apply the methods used in the bond percolation setting. One approach used in growing-degree expanders~\cite{random subgraphs}, similar to that used in the classical paper of Erd\H{o}s and R\'enyi~\cite{Erdos1} on the binomial random graph $G(n,p)$ and its close analogue $G(n,m)$ (see~\cite{JLR},~\cite{book random graphs} and~\cite{Intro to Random Graphs} for a systematic coverage), utilises that the edge-boundary of every (relatively small) subset of $G$ is about $d$ times the size of the subset. An analogous strategy in our settings requires that the neighbourhood of every subset inside $R$ is about $d$ times the size of the subset. However, as can be seen in the discussion following the statement of Lemma 2.4 here, this is not necessarily the case in site percolation on $(n,d,\lambda)$-graphs even for sublinear sets. In fact, the main technical part of our proof is to prove that this holds for small linear (in $n/d$) sets, and as it turns out, this suffices. Another approach, used in the constant-degree high-girth expanders~\cite{high-girth} and apparently originated in the classical paper of Ajtai, Koml\'os and Szemer\'edi~\cite{AKS}, would require us to establish the survival probability of every vertex (not an easy task in our setting, requiring a careful treatment of cycles) and then to continue with basic expansion of sets. In both approaches, the concentration results are obtained with the classical Azuma-Hoeffding bound --- which is too weak in our settings (see Section 5 of this paper). As such, we will require a careful combination of several approaches to obtain our results. A detailed comparison of our results with the analogous results in the $G(n,p)$ and bond percolation models follows each theorem.
	
	For a given set of vertices $V$, let $V_p$ be the random subset $V_p\subseteq V$ formed by including each vertex $v\in V$ independently and with probability $p$. We will make use of this notation throughout the rest of the paper.

	With this in mind, we can now state our main results. First, we improve the result of~\cite{site-percolation} for the subcritical phase. Notice that the statement does not assume anything on the structure of the $d$-regular graph, and is valid for any $d$-regular graph.
	\paragraph{Theorem 1} \emph{For every small enough constant $\epsilon>0$ and $d=o(n), d\ge 3$, let $p=\frac{1-\epsilon}{d}$. Let $G=(V,E)$ be a $d$-regular graph on $n$ vertices. Then, \textbf{whp}, all connected components in the induced subgraph $G[V_p]$ are of size at most}
	$$\frac{4}{\epsilon^2}\ln\left(\frac{n}{d}\right).$$

	For the size of a largest component in the supercritical phase, we note that the results we will present depend on the spectral ratio $\lambda/d$, which we will bound by $\delta$. As such, before presenting the theorems, we define the following notation:
	\begin{align}
	\delta(\alpha)=\alpha^{\frac{2}{\alpha}}.
	\end{align}
	Our results will involve accuracy parameter $\alpha$. In order to relate the results to the spectral ratio, we will assume $\frac{\lambda}{d}\le \delta(\alpha)=\alpha^{\frac{2}{\alpha}}$. We will make some mild assumptions on $\alpha$ with respect to $\epsilon$ (from $p=\frac{1+\epsilon}{d}$), and since the smallest spectral ratio is larger than $\frac{1}{\sqrt{d}}$ (see, for example,~\cite{expanders}), we may also assume that $\alpha=\omega\left(\frac{1}{\ln d}\right)$. In general, we allow $\alpha=\alpha(d)$ to be a function vanishing with $d\to\infty$.\\
	The extinction probability in an infinite $(d-1)$-ary tree rooted at a vertex $v$ in the site percolation model with probability $\frac{1+\epsilon}{d}$ satisfies:
	$$q=1-\frac{1+\epsilon}{d}+\frac{1+\epsilon}{d}q^{d-1},$$
	representing the probability of the vertex itself not being chosen, or that it was chosen yet all of its $d-1$ neighbours become extinct. We can then write $q=1-\frac{x}{d}$, and, using the approximation $(1-x/d)^d\approx\exp(-x)$ for large enough $d$, arrive at the following equation: $x=(1+\epsilon)\left(1-\exp\left(-x\right)\right)$.
	This equation has a unique solution in $(0,1)$. Indeed, the function $f(x)=x-(1+\epsilon)\left(1-\exp\left(-x\right)\right)$ is decreasing in the interval $[0,\ln(1+\epsilon)]$, and increasing afterwards; also, at $x=0$ it takes the value $0$, and at $x=1$ it is positive. Therefore, there is a unique $x$ in the interval $(\ln(1+\epsilon),1)$ where the function takes the value $0$. We denote by $x$ the unique solution in $(0,1)$ of:
	\begin{align}
	x=(1+\epsilon)\left(1-\exp\left(-x\right)\right).
	\end{align}
	Note that $x=2\epsilon-\frac{2\epsilon^2}{3}+O(\epsilon^3)$.
	
	We can now state the main result concerning the typical order of a largest component in the supercritical regime:
	\paragraph{Theorem 2} \emph{For every small enough constant $\epsilon>0$ and $d=o(n), d\ge 3$, let $p=\frac{1+\epsilon}{d}$. Let $2\sqrt{\frac{d}{n}}<\alpha<\epsilon^2$, and let $G=(V,E)$ be an $(n,d,\lambda)$-graph where $\lambda/d\le\delta(\alpha)$, with $\delta(\alpha)$ as defined in $(1)$. Denote by $L_1$ a largest connected component in $G[V_p]$. Then, there exists a positive constant $\beta=\beta(\epsilon)$ such that with probability at least $1-\exp\left(-\beta\frac{\alpha^2 n}{d}\right)$, 
	$$\Bigg||L_1|-\frac{xn}{d}\Bigg|\le\frac{7\alpha n}{d},$$
	where $x$ is from $(2)$.\\}
	
	Note that since $d/\lambda=O(\sqrt{d})$ for $d=o(n)$, our requirement $\lambda/d\le \delta(\alpha)$ means that the graph has at least large constant-degree, specifically $d=\Omega\left(\frac{1}{\delta^2}\right)$.

	In Theorem 3, we will also show that all the components other than this largest connected component are much smaller. Thus, Theorem 2 shows the size of \textit{the} largest component is concentrated around $\frac{xn}{d}$, and the tightness of the result depends on how small $\delta$ is, i.e. how strict is our assumption on the spectral ratio. To put this result in the right context, denote by $y$ the unique solution in $(0,1)$ of the equation
	\begin{align}
	y\exp\left(-y\right)=(1+\epsilon)\exp\left(-(1+\epsilon)\right),
	\end{align}
	(this definition was already present in the classical paper of Erd\H{o}s and R\'enyi~\cite{Erdos1}). In the $G(n,p)$ model and both bond percolation on constant-degree high-girth expanders and bond percolation on growing-degree expanders (with $p=\frac{1+\epsilon}{n}$ and $\frac{1+\epsilon}{d-1}$, respectively), the largest component is typically of size $\left(1-\frac{y}{1+\epsilon}\right)n$ (see Theorem 2.19 in~\cite{Intro to Random Graphs}, Theorem 1 in~\cite{high-girth} and Theorem 1 in~\cite{random subgraphs} respectively). In the site percolation model, we need to factor in the probability of a vertex to be included in $G[V_p]$ before looking at its progeny --- and indeed, rearranging $(3)$, we have that $y=(1+\epsilon)\exp(-1-\epsilon+y)$, that is $1+\epsilon-y=(1+\epsilon)\left(1-\exp(-1-\epsilon+y)\right)$, which means by $(2)$ that 
	\begin{align}
	    x=1+\epsilon-y.
	\end{align}
	Thus, the typical size of the largest connected component in the site percolation model is $\frac{1+\epsilon}{d}$ times the typical size of the largest connected component in the bond percolation model, which is what we would expect by this intuition.

	We also prove the following result for the typical size of the remaining components:
	\paragraph{Theorem 3} \emph{For every small enough constant $\epsilon>0$ and $d=o(n), d\ge 3$, let $p=\frac{1+\epsilon}{d}$. Assume that $\frac{2}{\ln(n/d)}<\alpha<\epsilon^4$. Let $G=(V,E)$ be an $(n,d,\lambda)$-graph where $\lambda/d\le\delta(\alpha)$, with $\delta(\alpha)$ as defined in $(1)$. Then, there exists a positive constant $\beta=\beta(\epsilon)$ such that with probability at least $1-\left(\frac{n}{d}\right)^{-\beta\ln\ln(n/d)}$, all but at most $\frac{15\alpha n}{d}$ of the vertices in $V_p$ belong either to isolated trees of order at most $\frac{1}{\alpha}$, or to the unique giant connected component in $G[V_p]$.}\\

	To put this result in the right context, recall that the largest connected component in the supercritical regime in both the $G(n,p)$ model and the bond percolation on growing-degree expanders model is known to be typically unique, while the remaining components are typically of size $O(\ln n)$ and are comprised mainly of fixed size trees. On the other hand, in the model of bond percolation on constant-degree high-girth expanders it was shown (see~\cite{high-girth}) that the second-largest component can be typically of size $n^\omega$ for $\omega$ any constant smaller than $1$ (this is in fact optimal, as can be seen by the results in the paper of Alon, Benjamini and Stacey~\cite{ABS}).  
	
	We also discuss some properties of the giant component:
	\paragraph{Theorem 4} \emph{For every small enough constant $\epsilon>0$ and $d=o(n), d\ge 3$, let $p=\frac{1+\epsilon}{d}$. Assume that $\frac{2}{\ln(n/d)}<\alpha<\epsilon^8$. Let $G=(V,E)$ be an $(n,d,\lambda)$-graph where $\lambda/d\le\delta(\alpha)$, with $\delta(\alpha)$ as defined in $(1)$. Let $e(L_1)$ denote the number of edges in the giant component in $G[V_p]$. Then, there exists a positive constant $\beta=\beta(\epsilon)$ such that with probability at least $1-\left(\frac{n}{d}\right)^{-\beta\ln\ln(n/d)}$, 
	$$\Big|e(L_1)-\frac{\left((1+\epsilon)^2-(1+\epsilon-x)^2\right)n}{2d}\Big|\le\frac{8\alpha^{1/4}n}{d},$$
	where $x$ is as defined in $(2)$.\\}

	Intuitively, for an edge to belong to the largest component, we require that both of its endpoints fall into $V_p$, and at least one of them develops a large component. Both endpoints fall into $V_p$ with probability $p^2$, and given that they are in $V_p$ and assuming independence, the probability that neither endpoint develops a large component is about $$\left(\frac{\frac{1+\epsilon-x}{d}}{\frac{1+\epsilon}{d}}\right)^2=\left(\frac{1+\epsilon-x}{1+\epsilon}\right)^2.$$
	Since $G$ is a $d$-regular graph and thus has $\frac{nd}{2}$ edges, we would anticipate the number of edges in the largest component to be about $$\frac{nd}{2}\left(\frac{1+\epsilon}{d}\right)^2\left(1-\left(\frac{1+\epsilon-x}{1+\epsilon}\right)^2\right)=\frac{\left((1+\epsilon)^2-(1+\epsilon-x)^2\right)n}{2d},$$
	and thus the result of Theorem 4 matches the intuition.
	
	The next theorem concerns the appearance of long cycles:
	\paragraph{Theorem 5} \textit{For every small enough constant $\epsilon>0$ and $d=o(n), d\ge 3$, let $p=\frac{1+\epsilon}{d}$. Assume that $2\sqrt{\frac{d}{n}}<\alpha<\epsilon^3$. Let $G=(V,E)$ be an $(n,d,\lambda)$-graph where $\lambda/d\le\delta(\alpha)$, with $\delta(\alpha)$ as defined in $(1)$. Then, there exists a positive constant $\beta=\beta(\epsilon)$ such that with probability at least $1-\exp\left(-\beta\frac{\alpha^2 n}{d}\right)$, there is a cycle of length at least $\frac{\epsilon^2n}{10^2d}$ in $G[V_p]$.\\}

	Note that for an edge to belong to a long cycle in $G[V_p]$, we require that both its endpoints fall into $V_p$ and that each one of them survives into a long path, not including this edge. This happens with probability about $\left(\frac{x}{d}\right)^2=\frac{\Theta(\epsilon^2)}{d^2}$. Considering all $nd/2$ edges in $G$, we can anticipate that a longest cycle in $G[V_p]$ will have $O(\epsilon^2)n/d$ edges.
	
	Finally, we treat typical expansion properties of the largest component. Here (and throughout the paper), we denote by $N_G(S)$ the external neighbourhood of the set $S$ in the graph $G$.
	\paragraph{Theorem 6} \textit{For every small enough constant $\epsilon>0$ and $d=o(n), d\ge 3$, let $p=\frac{1+\epsilon}{d}$. Let $2\sqrt{\frac{d}{n}}<\alpha<\epsilon^2$ and let $G=(V,E)$ be an $(n,d,\lambda)$-graph where $\lambda/d\le\delta(\alpha)$, with $\delta(\alpha)$ as defined in $(1)$. Then, there exists a positive constant $\beta=\beta(\epsilon)$ such that \textbf{whp} for every set of vertices $S$ belonging to the largest component in $G[V_p]$ with $\frac{16\alpha n}{d}\le |S| \le \frac{(x-9\alpha)n}{d}$, we have that $$\big|N_{G[V_p]}(S)\big|\ge\frac{\beta\alpha^2}{\ln\left(\frac{1}{\alpha}\right)}\cdot\frac{n}{d}.$$\\}

	Theorem 6 shows that the largest component is \textbf{whp} a reasonably good expander on the family of linearly-sized subsets. This can be used to derive other typical properties of the largest component --- we refer the reader to~\cite{expanders} for a survey including many results of this type.
	
	Our notation is fairly standard. We omit rounding signs for the sake of clarity of presentation.

\section{Auxiliary Lemmas}
\subsection{The DFS Algorithm}
	As it is crucial to our proofs, we will briefly discuss the DFS algorithm on a random vertex subgraph of a given $(n,d,\lambda)$-graph. Since the algorithm is well known and was discussed in the papers~\cite{K and S} and~\cite{site-percolation}, we will only briefly describe it. We define the following sets:
	\begin{description}
	\item[$\bullet$ $S$] is the set of vertices whose exploration is complete;
	\item[$\bullet$ $T$] is the set of unvisited vertices;
	\item[$\bullet$ $U$] is the set of currently explored vertices, kept in a stack;
	\item[$\bullet$ $W$] is the set of vertices discovered who fall outside of the random set $V_p$.
	\end{description}
	For a $G=(V,E)$ graph, the algorithm starts with $S=U=W=\emptyset$ and $T=V$, and ends when $U\cup T=\emptyset$. At each step, if $U$ is non-empty, the algorithm queries $T$ for neighbours of the last vertex in $U$, scanning these neighbours according to some 		$\sigma$ prioritization on the set $V$. If the last vertex in $U$ has a neighbour in $T$, the algorithm flips a coin with probability $p$. If the result of this coin flipping is positive, the algorithm moves the neighbour to $U$; otherwise, it moves it to $W$. If the last vertex in 		$U$ has no more neighbours in $T$, it moves to $S$. Finally, if $U$ is empty, the algorithm chooses a vertex from $T$ according to $\sigma$, and flips the coin to decide whether it moves to $U$ (positive) or to $W$ (negative). We feed the DFS algorithm with a 			sequence of i.i.d. Bernoulli$(p)$ random variables, $(X_i)_{i=1}^n$, so that the $i$-th coin flipping is answered positively if $X_i=1$, and negatively otherwise.\\
	Notice that the final subset $S$ of the algorithm is distributed exactly like a random subset $V_p$, formed by including each vertex of $V$ independently and with probability $p$. Furthermore, observe that at any stage of the algorithm, $S$ and $T$ have no edges between 		them, and therefore $N_G(S)\subseteq U\cup W$. Last, but not least, each connected component of $G[V_p]$ corresponds to an \emph{epoch} in the DFS run --- each epoch starts at the moment the first vertex enters $U$ and ends at the first subsequent moment where 		$U$ is empty once again.\\
	We note that in the application of the DFS algorithm to the case of bond percolation (see \cite{K and S}), there is a random variable corresponding to the number of queries between the stack $U$ and the set $T$, corresponding to edges whose coin flip was answered in the negative. Therefore, in order to obtain the asymptotic order of the giant in bond percolation utilising the DFS algorithm, one needs to estimate this random variable (see \cite{D and K} for an estimation of this random variable and a careful yet relatively simple analysis of the performance of the DFS in $G(n,p)$). Here, on the other hand, all the coin flips answered in the negative correspond to a vertex moving to $W$. This property will allow us to utilise this variant of the DFS to obtain the asymptotic order of the giant in this model in a more direct manner.

\subsection{Concentration of Random Variables}
	The DFS algorithm allows us to study the induced subgraph $G[V_p]$ via properties of the random sequence $(X_i)_{i=1}^n$, where we are specifically interested in epochs in the DFS run, as they correspond to connected components. Throughout the proofs of the 	theorems, we will require the following probabilistic lemma:
	\paragraph{Lemma 2.1} \emph{Let $\epsilon>0$ be a small enough constant and $d=o(n), d\ge 3$. Let $(X_i)_{i=1}^n$ be a sequence of i.i.d. Bernoulli$(p)$ random variables. Then, there exists a positive constant $\beta=\beta(\epsilon)$ such that the following is true:}
	\begin{itemize}
	\item[\emph{1.}] \emph{Let $p\le\frac{1+\epsilon}{d}$. Then, with probability at least $1-\exp\left(-\beta\frac{n}{d}\right)$ there are at most $\frac{2n}{d}$ random variables $X_i$ that took value 1.}
	\item[\emph{2.}] \emph{Let $p=\frac{1-\epsilon}{d}$ and $k=\frac{4}{\epsilon^2}\cdot\ln\left(\frac{n}{d}\right)$. Then, \textbf{whp} there is no interval $I$ of length of length $kd$ starting with a random variable that took value $1$ and in which at least $k$ of the random variables $X_i$, $i\in I$, took value 1.}
	\item[\emph{3.}] \emph{Let $p=\frac{1+\epsilon}{d}$ and let $c>0$. Then, with probability at least $1-\exp\left(-\beta\frac{c^2 n}{d}\right)$, for all $0\le t\le n$ such that $X_{t+1}=1$ the following holds:}
	$$\Bigg|\sum_{i=1}^{t}X_i-\frac{(1+\epsilon)t}{d}\Bigg|\le \frac{\epsilon^2 c n}{d}.$$
	\end{itemize}
	\begin{proof} \mbox{}
	\begin{itemize}
	\item[\emph{1.}] This is equivalent to stating that $\sum_{i=1}^n X_i\le \frac{2n}{d}$. This sum is distributed binomially with parameters $n$ and $p$. Using a Chernoff-type bound (see for example Theorem A.1.11 of~\cite{Chernoff from K and S}), we have that
	$$P\Bigg[Bin(n,p)\ge\frac{2n}{d}\Bigg]< \exp\left(-\frac{n}{7d}\right).$$
	\item[\emph{2.}] The sum $\sum_{i\in I}X_i$ is distributed binomially with parameters $kd$ and $p$. We denote the first random variable of the interval by $X_{i_1}$. Applying a Chernoff-type bound together with union bound over the $n$ possible different intervals, we have that the probability of an interval violating the assertion of the Lemma is at most:
	\begin{align*}
	n\cdot P[X_{i_1}=1]\cdot P\Bigg[Bin(kd-1,p)\ge k-1\Bigg]&<\frac{n}{d}\cdot\exp\left(-\frac{\epsilon^2(1-\epsilon)k}{3}\right)\\
	&<\frac{n}{d}\cdot\exp\left(-\frac{5}{4}\cdot\ln\left(\frac{n}{d}\right)\right)\\
	&=o(1),
	\end{align*}
	where the constant $\frac{1}{3}$ in the exponent of the first inequality comes from the Chernoff-type bound.
	\item[\emph{3.}] The proof is identical for both the lower tail and the upper tail, and we will thus show it for the lower tail. Applying a Chernoff-type bound for the lower tail of $Bin(t,p)$ together with union bound over the $n$ possible different values of $t$, we have that the probability for any $t$ violating the assertion of the Lemma is at most:
	\begin{align*}
	n\cdot P[X_{t+1}=1]\cdot P\Bigg[Bin(t,p)\le\frac{\left(1+\epsilon-\frac{\epsilon^2 c n}{t}\right)t}{d}\Bigg]&<\frac{2n}{d}\cdot\exp\left(-\frac{\epsilon^4 c^2 n}{3d}\right)\\
	&\le\exp\left(-\beta\frac{c^2 n}{d}\right).
	\end{align*}
	\end{itemize}
	\end{proof}

\subsection{Properties of Pseudo-Random Graphs}
    	In~\cite{site-percolation}, a key ingredient in the proofs for the lower bounds on the size of a longest path and a largest component in the supercritical phase was Lemma 3.1 therein, which provided a lower bound on the expansion of large enough sets. We will require tighter upper and lower bounds for the expansion of large enough sets. The idea behind the proof is very much alike to that of Lemma 3.1 of~\cite{site-percolation}.\\
	First, we will state the expander mixing lemma (see, for example, a somewhat stronger result of Theorem 2.11 of~\cite{survey pseudo}), followed by a short corollary. 
	\paragraph{Lemma 2.2} \emph{Let $G=(V,E)$ be an $(n,d,\lambda)$-graph. Then for every two subsets $B,C\subseteq V$,}
	$$\left|e(B,C)-\frac{d}{n}|B||C|\right|\le\lambda\sqrt{|B||C|}.$$
	\paragraph{Corollary 2.3} \emph{Let $G=(V,E)$ be an $(n,d,\lambda)$-graph, and let $\alpha>0$, $B\subseteq V$, $|B|\ge\frac{n}{2}$.
	\begin{itemize}
	\item[1.] Define:
	$$C=\left\{v\in V:d(v,B)\ge(1+\alpha)\frac{|B|d}{n}\right\},$$
	then $|C|\le \frac{2}{\alpha^2}\left(\frac{\lambda}{d}\right)^2n$.
	\item[2.] Define:
	$$C=\left\{v\in V:d(v,B)\le(1-\alpha)\frac{|B|d}{n}\right\},$$
	then $|C|\le\frac{2}{\alpha^2}\left(\frac{\lambda}{d}\right)^2n$, as well.
	\end{itemize}}
	\begin{proof}
	Note that Part 2 is exactly Corollary 2.2 in~\cite{site-percolation}, and therefore we will only prove Part 1. \\
	By the definition of $C$, we have $e(B,C)\ge(1+\alpha)\frac{|B||C|d}{n}$. On the other hand, by the expander mixing lemma (Lemma 2.2) we have $e(B,C)\le\frac{d}{n}|B||C|+\lambda\sqrt{|B||C|}$. Combining these we have:
	$$(1+\alpha)\frac{|B||C|d}{n}\le\frac{|B||C|d}{n}+\lambda\sqrt{|B||C|}.$$
	Recalling the fact $|B|\ge \frac{n}{2}$,
	$$|C|\le\frac{2}{\alpha^2}\left(\frac{\lambda}{d}\right)^2n,$$
	as required. 
	\end{proof}
	We can now state the key lemma, bounding the expansion of sets of relevant size.
	\paragraph{Lemma 2.4} \emph{Let $d=o(n), d\ge 3$, $G=(V,E)$ be an $(n,d,\lambda)$ graph. Let $2\sqrt{\frac{d}{n}}<\alpha<\epsilon^2$ and assume $\lambda/d\le\delta(\alpha)$, with $\delta(\alpha)$ as defined in $(1)$. Let $p\le2/d$. Then, there exists a positive constant $\beta$ such that with probability at least $1-\exp\left(-\beta\frac{\alpha n}{d}\right)$, $V_p$ does not contain a set $S$ of size $m$, $\frac{\alpha n}{d}\le m\le \frac{n}{3d}$, with:
	\begin{itemize}
	\item[1.]
	$$|N_G(S)|>(1+2\alpha)n\left(1-\exp\left(-\frac{dm}{n}\right)\right),$$
	or;
	\item[2.]
	$$|N_G(S)|<(1-2\alpha)n\left(1-\exp\left(-\frac{dm}{n}\right)\right).\\$$
	\end{itemize}}
	
	Before we prove Lemma 2.4, we should mention that we indeed cannot bound in a similar manner the expansion of sublinear-sized subsets of $V(G)$. To see that, let $G_0$ be an $(n_0,d_0,\lambda_0)$-graph with $\lambda_0/d_0\le\delta$ and $\lambda_0=\Theta(\sqrt{d_0})$, and therefore $d_0=\Theta\left(\frac{1}{\delta^2}\right)$. Consider the blow-up graph $G=G_0(2)$, where we replace each vertex of $G_0$ by an independent set of size $2$ and connect two vertices of $G$ by an edge if and only if the corresponding vertices of $G_0$ are connected by an edge. $G$ is then an $(n,d,\lambda)$-graph with $n=2n_0$, $d=2d_0$ and $\lambda=2\lambda_0$ (see for example Proposition 2.5 of~\cite{triangle expanders}). 
	Denote the independent sets of $G$ (i.e. the images of the vertices of $G_0$) by $B_1$ up to $B_{n_0}$, and denote by $B(v)$ the independent set which contains the vertex $v$. Then,
	\begin{align*}
	    P\left[\big|B_i\cap V_p\big|= 2\right]& = p^2\\
	    &=\Theta\left(\frac{1}{d_0^2}\right)\\
	    &=\Theta\left(\delta^4\right).
	\end{align*}
	Therefore, if we define $S=\Big\{v\in V_p: \big|B(v)\cap V_p\big|= 2\Big\}$, then typically $|S|=\Theta\left(np\delta^4)=\Theta(\delta^4\frac{n}{d}\right)$. Furthermore, since we can couple the vertices of $S$ into pairs where every pair belongs to the same independent set, and therefore has the same set of $d$ neighbours, we have that $|N(S)|\le |S|d/2$.
	
	We are now ready to prove Lemma 2.4:
	
	\begin{proof}
	We will require the following inequalities:
	\begin{equation}
	1+x\le \exp(x), \hspace{2em} \forall x,
	\end{equation}
	\begin{equation}
	\hspace{4em}1-x\ge \exp\left(\frac{-x}{1-x}\right), \hspace{2em} 0\le x<1,
	\end{equation}
	which are Lemma 22.1 in~\cite{Intro to Random Graphs}. We will also use the fact that:
	\begin{equation}
	\hspace{4em}\exp\left(-x\right)\le1-x+\frac{x^2}{2}, \hspace{2em} x\ge 0,
	\end{equation}
	which is a direct consequence of the representation of $\exp(x)$ as a power series. 
	\begin{itemize}
	\item[\emph{1.}] We call an $m$-set $S\subset V$ \emph{over-expanding} if $$|N_G(S)|>(1+2\alpha)n\left(1-\exp\left(-\frac{dm}{n}\right)\right).$$  \\
	We consider the number of ways to choose a sequence of distinct vertices of $G$, $\tau=(v_1,\cdots,v_m)$, such that the union $S$ of the vertices in the sequence forms an over-expanding set. Assume we have chosen the first $i-1$ vertices of $\tau$, and define $S_{i-1}=\{v_1,\cdots,v_{i-1}\}$ and $N_{i-1}=N_G(S_{i-1})$. We call a vertex $v$ \emph{bad} with respect to the prefix $(v_1,\cdots,v_{i-1})$ if $v$ has at least $(1+\alpha)\frac{d}{n}\left(n-|S_{i-1}\cup N_{i-1}|\right)$ neighbours in $V-(S_{i-1}\cup N_{i-1})$, and \emph{good} otherwise. Each good vertex $v_i$ appended to $S_{i-1}$ tames the increase in size of the external neighbourhood. We will now show that if $\tau$ has at most $\alpha m$ bad vertices, it cannot be over-expanding. \\
	We consider the worst case scenario, that is, where $N_G(S)$ has the largest number of vertices, given that $\alpha m$ are bad and $(1-\alpha)m$ are good.\\
	In this scenario, all the $\alpha m$ bad vertices are the last to join the sequence, since otherwise their drastic contribution to the increment of the size of the neighbourhood will constrain the possible increment to the size of the neighbourhood that a good vertex can contribute. In this scenario, all the $\alpha m$ bad vertices add, at most, $d$ new neighbours to $N$, that is $\alpha m d$ neighbours altogether at the end of the sequence. \\
	We will now bound the size of the neighbourhood of the $(1-\alpha)m$ good vertices, which are the first in the sequence. We have that $|S_0|,|N_0|=0$ and $|S_i|=|i|$, $|N_i|\le |N_{i-1}|+(1+\alpha)\left(d-\frac{d}{n}|N_{i-1}|\right)$. We will show by induction that $|N_k|\le n\left(1-\left(1-\frac{(1+\alpha)d}{n}\right)^k\right)$. Indeed, for $k=0$, this holds by definition. Assume it holds for some $k$, then we have:
	\begin{align*}
	|N_{k+1}|&\le|N_k|+(1+\alpha)\left(d-\frac{d}{n}|N_k|\right) \\
	&\le n\left(1-\left(1-\frac{(1+\alpha)d}{n}\right)^k\right)+(1+\alpha)\left(d-d\left(1-\left(1-\frac{(1+\alpha)d}{n}\right)^k\right)\right)\\
	&=n-n\left(1-\frac{(1+\alpha)d}{n}\right)^k+(1+\alpha)d\left(1-\frac{(1+\alpha)d}{n}\right)^k\\
	&=n\left(1-\left(1-\frac{(1+\alpha)d}{n}\right)\left(1-\frac{(1+\alpha)d}{n}\right)^k\right)\\
	&=n\left(1-\left(1-\frac{(1+\alpha)d}{n}\right)^{k+1}\right),
	\end{align*}
	where we used the fact that $d<\frac{n}{1+\alpha}$, and therefore taking the maximal value of $|N_k|$ leads to the maximal value of the expression. We thus have:
	$$|N_G(S)|\le \alpha dm + n\left(1-\left(1-\frac{(1+\alpha)d}{n}\right)^{(1-\alpha)m}\right).$$
	By inequality $(6)$:
	\begin{align*}
	\left(1-\frac{(1+\alpha)d}{n}\right)^{(1-\alpha)m}&\ge\exp\left(-(1-\alpha)m\left(\frac{\frac{(1+\alpha)d}{n}}{1-\frac{(1+\alpha)d}{n}}\right)\right)\\
								&=\exp\left(-\frac{\frac{(1-\alpha^2)dm}{n}}{1-\frac{(1+\alpha)d}{n}}\right)\\
						&\ge \exp\left(-\frac{dm}{n}\right),
	\end{align*}
	where the last inequality holds since $\frac{(1+\alpha)d}{n}<\alpha^2$, as we required $\alpha>2\sqrt{\frac{d}{n}}$. Furthermore, by inequality $(7)$: 
	\begin{align*}
	2\alpha n\left(1-\exp\left(-\frac{dm}{n}\right)\right)&\ge 2\alpha\left(dm-\frac{d^2m^2}{2n}\right)\\
	&\ge \alpha dm,
	\end{align*}
	where the last inequality holds since $m<\frac{n}{3d}$.
	Returning to $|N_G(S)|$, we now know that:
	\begin{align*}
	|N_G(S)|&\le \alpha dm + n\left(1-\left(1-\frac{(1+\alpha)d}{n}\right)^{(1-\alpha)m}\right)\\
		&\le \alpha dm + n\left(1-\exp\left(-\frac{dm}{n}\right)\right)\\
	    &\le (1+2\alpha)n\left(1-\exp\left(-\frac{dm}{n}\right)\right),
	\end{align*}
	that is, in order for $S$ to be an over-expanding set, the sequence $\tau$ must have at least $\alpha m$ bad vertices. Furthermore, we have for $i\le m$ that the set $S_{i-1}\cup N_{i-1}$ could have no more than $(i-1)(d+1)<m(d+1)<n/2$ vertices (since $m\le\frac{n}{3d}$).

	We can now use Corollary 2.3, and conclude that the number of bad choices for $v_i$ is at most $\frac{2}{\alpha^2}\delta^2n$. Therefore, the number of sequences $\tau$ 		with at least $\alpha m$ bad vertices is at most
	$${m\choose \alpha m}\left(\frac{2}{\alpha^2}\delta^2n\right)^{\alpha m}n^{m-\alpha m}\le \left(\left(\frac{e}{\alpha}\right)^\alpha\left(\frac{2}{\alpha^2}\delta^2\right)^\alpha n\right)^m.$$
	Dividing by $m!$ to get the number of unordered sets of size $m$ which could violate this property, and multiplying by $p^m$ to get the probability that $V_p$ contains a set of size $m$ violating this property, we have that this probability is at most
	$$\left(\left(\frac{2e\delta^2}{\alpha^3}\right)^\alpha \frac{enp}{m}\right)^m.$$
	Since we assume that $p\le \frac{2}{d}$ and $m\ge\frac{\alpha n}{d}$, it follows that $\frac{np}{m}\le\frac{2}{\alpha}$. We can now choose $\delta=\delta(\alpha)$ according to $(1)$, and the above probability will be at most $2^{-m}$. This, together with the union bound over all 		possible values of $m$ between $\frac{\alpha n}{d}$ and $\frac{n}{3d}$, completes the proof.

	\item[\emph{2.}] We call an $m$-set $S\subset V$ \emph{under-expanding} if $$|N_G(S)|<(1-2\alpha)n\left(1-\exp\left(-\frac{dm}{n}\right)\right).$$ \\
	Mirroring the proof of Part 1, we consider the number of ways to choose a sequence $\tau$ of $m$ distinct vertices of $G$, such that the union $S$ of the vertices in $\tau$ forms an under-expanding set. We define $S_{i-1}, N_{i-1}$ as in Part 1 above. 		We call a vertex $v$ \emph{bad} with respect to the prefix $(v_1,\cdots,v_{i-1})$ if $v$ has at most $(1-\alpha)\frac{d}{n}\left(n-|S_{i-1}\cup N_{i-1}|\right)$ neighbours in $V-(S_{i-1}\cup N_{i-1})$, and \emph{good} otherwise. This time, each good vertex $v_i$ 			appended to $S_i$ 	increases substantially the size of the external neighbourhood. Mirroring the case in Part 1, we will now show that if $\tau$ has at most $\alpha m$ bad vertices, it cannot be under-expanding. \\
	We consider the worst case scenario, which is when the set $N_G(S)$ has the smallest number of vertices, given that $\alpha m$ are bad and $(1-\alpha)m$ are good. \\
	All $\alpha m$ bad vertices should be the last to join the sequence, as otherwise they would allow for a larger increment in the external neighbourhood from the good vertices. In the worst case scenario, each bad vertex will add no new neighbours. As for the good 		vertices in the sequence, we have $|S_0|, |N_0|=0$, and $|S_i|=i$, $|N_i|\ge |N_{i-1}|+(1-\alpha)\frac{d}{n}\left(n-|S_{i-1}\cup N_{i-1}|\right)$. We will show by induction that $$|N_k|\ge n\left(1-\left(1-\frac{(1- 		\alpha)d}{n}\right)^k\right)-\frac{k(k-1)(1-\alpha)d}{2n}.$$ Indeed, for $k=0$, this holds by definition. Assume this holds for some $k$, then we have:
	\begin{align*}
	|N_{k+1}|&\ge|N_k|+(1-\alpha)\left(d-\frac{d}{n}|N_k\cup S_k|\right)\\
	&\ge n\left(1-\left(1-\frac{(1-\alpha)d}{n}\right)^k\right)-\frac{k(k-1)(1-\alpha)d}{2n}\\
	&\hspace{6em}+(1-\alpha)\left(d-\frac{d}{n}\left(k+n\left(1-\left(1-\frac{(1-\alpha)d}{n}\right)^k\right)\right)\right)\\
	&=n-\left(n-(1-\alpha)d\right)\left(1-\frac{(1-\alpha)d}{n}\right)^k-\frac{(1-\alpha)d}{n}\left(\frac{k(k-1)}{2}+k\right)\\
	&=n\left(1-\left(1-\frac{(1-\alpha)d}{n}\right)^{k+1}\right)-\frac{(k+1)k(1-\alpha)d}{2n},
	\end{align*}
	where we used the fact that $d<n$, and therefore taking the minimal value of $|N_k|$ leads to the minimal value of the expression. We thus have:
	$$|N_G(S)|\ge n\left(1-\left(1-\frac{(1-\alpha)d}{n}\right)^{(1-\alpha)m}\right)-\frac{(1-\alpha)^3dm^2}{2n}.$$
	By inequality $(5)$:
	\begin{align*}
	\left(1-\frac{(1-\alpha)d}{n}\right)^{(1-\alpha)m}\le\exp\left(-\frac{(1-\alpha)^2dm}{n}\right),
	\end{align*}
	and therefore,
    $$|N_G(S)|\ge n\left(1-\exp\left(-\frac{(1-\alpha)^2dm}{n}\right)\right)-\frac{(1-\alpha)^3dm^2}{2n}.$$
    Now, consider the function:
    \begin{align*}
    f(m)&=n\left(1-\exp\left(-\frac{(1-\alpha)^2dm}{n}\right)\right)-\frac{(1-\alpha)^3dm^2}{2n}\\
    &\hspace{18em}-(1-2\alpha)n\left(1-\exp\left(-\frac{dm}{n}\right)\right)\\
    &=2\alpha n+(1-2\alpha)n\exp\left(-\frac{dm}{n}\right)-n\exp\left(-\frac{(1-\alpha)^2dm}{n}\right)-\frac{(1-\alpha)^3dm^2}{2n}.
    \end{align*}
    Observing its derivative, we note that $f(m)$ increases with $m$, and therefore:
    \begin{align*}
    f(m)&\ge f(0)\\
        &=-2\alpha n+(1-2\alpha)n-n\\
        &=0.
    \end{align*}
    To conclude:
    \begin{align*}
    |N_G(S)|&\ge n\left(1-\exp\left(-\frac{(1-\alpha)^2dm}{n}\right)\right)-\frac{(1-\alpha)^3dm^2}{2n}\\
			&\ge(1-2\alpha)n\left(1-\exp\left(-\frac{dm}{n}\right)\right).
    \end{align*}
	As in Part 1, here too we have that the set $S_{i-1}\cup N_{i-1}$ could have no more than $n/2$ vertices, and the bad vertices were defined in accordance with Corollary 2.3. To complete the proof, we simply repeat the computation for the number of $m$-sets violating the Lemma, where the computation is the same as in Part 1.
	\end{itemize}
	\end{proof}

\section{The size of Components in the Subcritical Phase}
	The proof is similar to that in~\cite{site-percolation}.
	\begin{proof}[\emph{\textbf{Proof of Theorem 1}}]
	Assume to the contrary that $G[V_p]$ contains a component of size at least $k=\frac{4}{\epsilon^2}\ln\left(\frac{n}{d}\right)$. Let us observe the epoch in the DFS where this component was discovered. There is a moment in this epoch, where the algorithm 		found the $k$-th vertex of the component and has just moved it into $U$. Denote by $C_0$ the portion of the component discovered by that moment, i.e. all the vertices of this component which are already in $S\cup U$. Then $|C_0|=k$, and the subgraph $G[C_0]$ is 		connected and spans at least $k-1$ edges. We thus have that:
	$$|N_G(C_0)|\le e_G(C_0, V\backslash C_0)\le kd-2(k-1).$$
	Observe that exactly $k$ random variables $X_i$ took value $1$ during the epoch at that moment, and only the vertices in $C_0$ and those neighbouring them in $G$ have been queried. That means that we have had at most $kd+k-2(k-1)\le kd$ queries from the beginning of the epoch, and in them we have had at least $k$ random variables $X_i$ that took value 1. Furthermore, each epoch starts at some $X_{i_1}$ with $X_{i_1}=1$. Thus, we may conclude that this is a contradiction to Property 2 of Lemma 2.1.
	\end{proof}

\section{The Size of the Largest Component in the Supercritical Phase}
	We begin with an upper bound. 
	\begin{proof}[\emph{\textbf{Proof of Theorem 2 --- upper bound}}]
	Set $k=\left(x+6\alpha\right)n/d$. Assume to the contrary that $G[V_p]$ contains a component of size at least $k$, and consider the epoch in the DFS where this component was discovered. As in the proof of Theorem 1, there is a moment in this epoch where the $k$-th vertex of the component has just moved into $U$. At that moment, $k$ random variables $X_i$ took value $1$, while we we could only query these vertices and their neighbours in $G$. Invoking Property 1 of Lemma 2.4 with $\delta(\alpha)$ and $m:=k$ (note that $\frac{\alpha n}{d}\le k \le \frac{n}{3d}$), we have that with probability at least $1-\exp\left(-\beta'\frac{\alpha n}{d}\right)$ these $k$ random variables took value $1$ in an interval of length at most
	\begin{align*}
	k+(1+2\alpha)n\left(1-\exp\left(-\frac{dk}{n}\right)\right),
	\end{align*}
	and for our $k$ the above is
	$$k+(1+2\alpha)n\left(1-\exp\left(-x-6\alpha\right)\right).$$
	Observe that by the definition of $x$:
	\begin{align*}
	1-\exp\left(-x-6\alpha\right)&=\frac{x}{1+\epsilon}+\exp\left(-x\right)-\exp\left(-x-6\alpha\right)\\
	&=\frac{x}{1+\epsilon}+\exp(-x)\left(1-\exp\left(-6\alpha\right)\right)\\
	&=\frac{x}{1+\epsilon}+\left(1-\frac{x}{1+\epsilon}\right)\left(1-\exp(-6\alpha)\right).
	\end{align*}
	Recalling that according to $(2)$ we have that $x=2\epsilon-\frac{2\epsilon^2}{3}+O(\epsilon^3)$, we have for small enough $\epsilon$:
	\begin{align*}
	\frac{x}{1+\epsilon}+\left(1-\frac{x}{1+\epsilon}\right)\left(1-\exp(-6\alpha)\right)<\frac{x}{1+\epsilon}+6\alpha\left(1-\frac{7\epsilon}{4}\right).
	\end{align*}
	To conclude, the number of queries is at most:
	$$k+(1+2\alpha)\left(\frac{x}{1+\epsilon}+6\alpha\left(1-\frac{7\epsilon}{4}\right)\right)n.$$
	As in the proof of Theorem 1, we note that the positive answers to the queries in the interval are stochastically dominated by 
	$$Bin\left(k+(1+2\alpha)\left(\frac{x}{1+\epsilon}+6\alpha\left(1-\frac{7\epsilon}{4}\right)\right)n,\frac{1+\epsilon}{d}\right),$$
	and the expectation is at most (again, according to $(2)$): 
	\begin{align*}
	\mu&\le\frac{(x+3\alpha)n}{d}+\frac{(2\alpha x)n}{d}-\frac{9\alpha \epsilon n}{2d}\\
	&\le\frac{(x+3\alpha)n}{d}-\frac{\alpha\epsilon n}{2d}.
	\end{align*}
	Using a standard Chernoff-type bound, we have that the probability of such an event is at most:
	\begin{align*}
	P\Bigg[Bin\left(k+(1+2\alpha)\left(\frac{x}{1+\epsilon}+6\alpha\left(1-\frac{7\epsilon}{4}\right)\right)n,\frac{1+\epsilon}{d}\right)>k\Bigg]&<\exp\left(-\frac{\left(\frac{\alpha\epsilon n}{2d}\right)^2}{\frac{3(x+6\alpha)n}{d}}\right)\\
	&\le\exp\left(-\frac{\alpha^2\epsilon n}{25d}\right).
	\end{align*}
	By Property 1 of Lemma 2.1, there are, with probability $1-\exp\left(-\frac{n}{7d}\right)$, at most $\frac{2n}{d}$ such intervals (since each interval starts with a random variable that took value $1$). Using the union bound over these intervals and Property 1 of Lemma 2.1, we conclude that with probability at least $1-\exp\left(-\beta\frac{\alpha^2 n}{d}\right)$, there is no connected component of size at least $k$ in $G[V_p]$.
	\end{proof}

	We proceed with the lower bound, the proof of which is similar to that of Theorem 2 in~\cite{site-percolation}:
	\begin{proof}[\textbf{\emph{Proof of Theorem 2 --- lower bound}}] 
	Assume, for the sake of contradiction, that at some moment $$t\in\Big[\frac{3\alpha n}{2}, \frac{(x-5\alpha)n}{1+\epsilon}\Big],$$ we have that $U$ empties. We then have that $|S\cup W|=t$ and $$m:=|S|=\sum_{i=1}^tX_i\ge\frac{(1+\epsilon)t}{d}-\frac{\epsilon^2\alpha n}{d},$$ with probability at least $1-\exp\left(-\beta'\frac{\alpha^2 n}{d}\right)$ by Property 3 of Lemma 2.1 (using $\alpha$ instead of $c$). Since $\frac{\alpha n}{d}\le m\le\frac{n}{3d}$, we may invoke Lemma 2.4 with $\delta(\alpha)$. Recalling that $U=\emptyset$ and thus $N_G(S)\subseteq W$, we have by property 2 of Lemma 2.4 that:
	\begin{align*}
	|W|&\ge(1-2\alpha)n\left(1-\exp\left(-\frac{m\cdot d}{n}\right)\right)\\
	&\ge (1-2\alpha)n\left(1-\exp\left(-\frac{(1+\epsilon)t}{n}+3\epsilon^2\alpha \right)\right).
	\end{align*}
	Since $|W|<t$, it suffices to show that $(1-2\alpha)n\left(1-\exp\left(-\frac{(1+\epsilon)t}{n}+3\epsilon^2\alpha \right)\right)>t$ in order to obtain the contradiction. Note that $$f(t)=(1-2\alpha)n\left(1-\exp\left(-\frac{(1+\epsilon)t}{n}+3\epsilon^2\alpha \right)\right)-t$$ is a concave function, and as such it is sufficient to prove that $f(t)>0$ for the left and right ends of the interval. For $t=\frac{3\alpha n}{2}$, using inequality $(7)$, we obtain:
	\begin{align*}
	    f\left(\frac{3\alpha n}{2}\right)&\ge(1-2\alpha)n\left(1-\exp\left(-\frac{3(1+\epsilon)\alpha}{2}+3\epsilon^2\alpha \right)\right)-\frac{3\alpha n}{2}\\
	    &\ge (1-2\alpha)n\left(1-\exp\left(-\frac{3\alpha}{2}-\epsilon\alpha\right)\right)-\frac{3\alpha n}{2}\\
	    &\ge(1-2\alpha)n\left(\frac{3\alpha}{2}+\epsilon\alpha-4\epsilon\alpha^2\right)-\frac{3\alpha n}{2}\\
	    &>0,
	\end{align*}
	where the last two inequalities hold since $\alpha<\epsilon^2$. As for $t=\frac{(x-5\alpha)n}{1+\epsilon}$, by the definition of $x$ and by the inequality $\exp(c)<1+c+c^2$ (for small enough $c>0$), we obtain:
	\begin{align*}
	    f\left(\frac{(x-5\alpha)n}{1+\epsilon}\right)&\ge
	    (1-2\alpha)n\left(1-\exp\left(-x+5\alpha+3\epsilon^2\alpha\right)\right)-\frac{(x-5\alpha)n}{1+\epsilon}\\
	    &\ge (1-2\alpha)n\left(\frac{x}{1+\epsilon}+\left(1-\frac{7\epsilon}{4}\right)\left(1-\exp\left(5\alpha+3\epsilon^2\alpha\right)\right)\right)-\frac{(x-5\alpha)n}{1+\epsilon}\\
	    &\ge(1-2\alpha)n\left(\frac{x}{1+\epsilon}+\left(1-\frac{7\epsilon}{4}\right)(-5\alpha-3\epsilon^2\alpha-10\alpha^2)\right)-\frac{(x-5\alpha)n}{1+\epsilon}\\
	    &\ge \frac{x-5\alpha}{1+\epsilon}-\frac{2\alpha x+3\epsilon\alpha}{1+\epsilon}+8\epsilon\alpha-\frac{(x-5\alpha)n}{1+\epsilon}\\
	    &>0,
	\end{align*}
	where we once again used the fact that $\alpha<\epsilon^2$. 
	We thus have that $U$ does not empty in the interval $\Big[\frac{3\alpha n}{2}, \frac{(x-5\alpha)n}{1+\epsilon}\Big]$, and all the positive answers between these two moments belong to the same connected component, whose size is at least:
	$$\sum_{i=1}^{\frac{(x-5\alpha)n}{1+\epsilon}}X_i-\sum_{i=1}^{\frac{3\alpha n}{2}}X_i\ge\frac{(x-7\alpha)n}{d},$$
	with probability at least $1-\exp\left(-\beta\frac{\alpha^2n}{d}\right)$, using a standard Chernoff-type bound.
	\end{proof}

\section{The Size of the Remaining Components}
    Throughout this and the next section, unless stated otherwise, we assume that $d=o(n), d\ge 3$ and $p=\frac{1+\epsilon}{d}$, and let $G=(V,E)$ be an $(n,d,\lambda)-$graph, with $\lambda/d\le\delta(\alpha)$ with $\delta(\alpha)$ as defined in $(1)$.

    We require the following Lemma, which holds for any $d$-regular graph:
    \paragraph{Lemma 5.1} \textit{The number of $k$-vertex trees contained in a $d$-regular graph on $n$ vertices is at least:
    $$n\frac{k^{k-2}(d-k)^{k-1}}{k!}.$$}
    This is Lemma 2 of~\cite{trees}.\\

    Note that Lemma 5.1 counts two trees with the same vertices, but different sets of edges, as two different trees. We want to bound from below the number of trees with different set of vertices. Let $G=(V,E)$ be a graph. We call a set of vertices $S\subseteq V(G)$ a \textit{connected }$k$-\textit{set} if $|S|=k$ and $G[S]$ is connected. We call such a connected $k$-set \textit{acyclic} if $G[S]$ contains no cycles.
    \paragraph{Lemma 5.2} \textit{Let $\alpha\ge \frac{2}{\ln\left(\frac{n}{d}\right)}$, $k\le \frac{1}{\alpha}$. Denote by $t_k$ the number of acyclic connected $k$-sets in $G$. Then:
    $$t_k\ge (1-\alpha)n\frac{k^{k-2}d^{k-1}}{k!}.$$}
    \begin{proof}
    Any connected $k$-set contains a spanning tree. By Lemma 5.1, we have at least $n\frac{k^{k-2}(d-k)^{k-1}}{k!}$ trees on $k$ vertices, with different sets of edges. Two such trees can have the same vertices only if the graph induced by their vertices contains a cycle. 
    
    Each cycle of length $\ell$ is composed of a path of $uv$ of length $\ell-3$, and two neighbours, one of $u$ and one of $v$, which are connected by an edge. Starting with a vertex $u$, we have at most $d^{\ell-3}$ paths of length $\ell-3$. In order to close such a path into a cycle, we consider the neighbourhood of $u$ and $v$ (disjoint from the path), each of size at most $d$, and then the number of edges between these neighbourhoods, which is by the expander mixing lemma (Lemma 2.2) at most $$\frac{d}{n}d^2+\lambda d=\frac{d^3}{n}+\lambda d.$$ 
    Considering all the $<d^{\ell-3}$ paths starting with a vertex $v$, all the $n$ different $v$'s we can start from, and the fact that we over-count each cycle $\ge \ell$ times (since we can start the cycle from $\ell$ different vertices), we conclude that we have at most:
    $$\frac{nd^{\ell-3}\left(\frac{d^3}{n}+\lambda d\right)}{\ell}=\frac{d^\ell+n\lambda d^{\ell-2}}{\ell},$$
    cycles of length $\ell$.
    
    Fix a cycle of length $\ell$, denote it by $C_\ell$. We want to bound the number of $k$-vertex trees whose vertex set contains $C_\ell$. For that, we follow the method  of Lemma 2 of~\cite{trees}. We fix a labelling of the vertices of $C_\ell$, $f_0:V(C_\ell)\to \{1,2,\cdots,\ell\}$. Note that for each $C_\ell$ we have $\ell!$ different labels as such. Given a tree $T$ on $k$ vertices whose vertex set contains $C_\ell$, we define $f:V(T)\to\{1,2,\cdots,k\}$ to be a labelling that extends $f_0$. Consider the pairs $(T,f)$ where $T$ is a $k$-vertex tree containing $C_\ell$, and $f$ is a labelling that extends $f_0$. Clearly, each $k$-vertex tree $T$ containing $C_\ell$ is in $(k-\ell)!$ such pairs. Furthermore, each such pair defines a unique labelled spanning tree $T'$ of $K_k$, where $(i,j)$ is an edge of $T'$ if and only if there is an edge $xy$ of $T$ such that $f(x)=i, f(y)=j$. Run a DFS on $T'$, starting at some vertex of $C_\ell$, and on reaching a vertex $m$ define $f^{-1}(m)$. For the $\ell$ vertices of $C_\ell$, this is predetermined. For each of the $k-\ell$ other vertices, there will be at most $d$ choices. Since there are $k^{k-2}$ labelled spanning trees of $K_k$, we can conclude that the number of $k$-vertex trees whose vertex-set contains $C_\ell$ is at most: $$\frac{d^{k-\ell}k^{k-2}}{(k-\ell)!\ell!},$$
    and thus, the number of such sets containing a cycle of length $\ell$, $3\le \ell\le k$, is at most:
    $$\frac{d^\ell+n\lambda d^{\ell-2}}{\ell}\cdot\frac{d^{k-\ell}k^{k-2}}{(k-\ell)!\ell!}=\frac{k^{k-2}}{\ell(k-\ell)!\ell!}\left(d^k+n\lambda d^{k-2}\right).$$

    Each $k$-vertex tree is either a unique acyclic connected $k$-set, or its corresponding connected $k$-set spans a cycle. We thus obtain the following bound:
    \begin{align*}
    t_k&\ge n\frac{k^{k-2}(d-k)^{k-1}}{k!}-\sum_{\ell=3}^{k}\frac{k^{k-2}}{\ell(k-\ell)!\ell!}\left(d^k+n\lambda d^{k-2}\right)\\
    &\ge n\frac{k^{k-2}(d-k)^{k-1}}{k!}-\frac{2^k k^{k-2}}{k!}\left(d^k+n\lambda d^{k-2}\right)\\
    &= n\frac{k^{k-2}}{k!}\left((d-k)^{k-1}-2^k\left(\frac{d^k}{n}+\lambda d^{k-2}\right)\right),
    \end{align*}
    where we used the identity $\sum_{x=0}^k\frac{1}{(k-x)!x!}=\frac{2^k}{k!}$. Note that by Bernoulli's inequality, $$(d-k)^{k-1}>d^{k-1}\left(1-\frac{k}{d}\right)^k\ge d^{k-1}\left(1-\frac{k^2}{d}\right),$$
    and we can conclude:
    \begin{align*}
    t_k&\ge n\frac{k^{k-2}d^{k-1}}{k!}\left(1-\frac{k^2}{d}-2^k\left(\frac{d}{n}+\frac{\lambda}{d}\right)\right)\\
    &\ge (1-\alpha)n\frac{k^{k-2}d^{k-1}}{k!},
    \end{align*}
    where the last inequality is due to our assumptions on $\alpha$, $k$ and by $(1)$.
    \end{proof}

    We want to show now that the number of acyclic connected $k$-sets in $G[V_p]$ is tightly concentrated. For that, we will use of the following (stronger) variant of the Azuma-Hoeffding inequality due to Warnke~\cite{Stronger Azuma}:
    \paragraph{Lemma 5.3} \textit{Let $1\le i\le N$. Let $X=(X_1,\cdots,X_N)$ be a family of random variables with $X_i$ taking values in a set $\Delta_i$. Let $\Gamma\subseteq \Pi_{j\in[N]}\Delta_j$ be an event. Assume that the function ${f:\Pi_{j\in[N]}\Delta_j\to \mathbb{R}}$ satisfies that there are numbers $(c_i)_{i\in[N]}$ and $(d_i)_{i\in [N]}$ with $c_i\le d_i$ such that the following holds for any two possible sequences of outcomes $a_1,\cdots,a_{i-1},a$ and $a_1,\cdots,a_{i-1},b$ of $X_1,\cdots,X_i$. Defining for $z\in \Delta_i$
    $$\Sigma_z=\left\{x=(a_1,\cdots,a_{i-1},z,x_{i+1},\cdots,x_N)\in\Pi_{j\in[N]}\Delta_j:P[X=x]>0\right\},$$
    there is an injection $\rho_i=\rho_i\left(\Sigma_a,\Sigma_b\right):\Sigma_a\to\Sigma_b$ such that for all $x\in\Sigma_a$ we have 
    \begin{align*}
        &\Big|f(x)-f\left(\rho_i(x)\right)\Big|\le 
        \begin{cases}
        c_i & \text{if } x\in \Gamma, \\
        d_i & \text{otherwise,}
        \end{cases} \text{ and,}\\
        &P\left[X=x|X\in\Sigma_a\right]\le P\left[X=\rho_i(x)|X\in\Sigma_b\right].
    \end{align*}
    Then, for any numbers $(\gamma_i)_{i\in [N]}$ with $\gamma_i\in (0,1]$, there is an event $\mathcal{B}=\mathcal{B}\left(\Gamma,(\gamma_i)_{i\in[N]}\right)$ satisfying
    $$P[\mathcal{B}]\le \sum_{i=1}^N\gamma_i^{-1}P[X\notin \Gamma] \hspace{2em} \text{and} \hspace{2em} \neg\mathcal{B}\subseteq \Gamma,$$
    such that for $e_i=\gamma_i(d_i-c_i)$, and any $t\ge 0$ we have:
    $$P[f(x)\le \textbf{E}f(x)-t \text{ and } \neg\mathcal{B}]\le \exp\left(-\frac{t^2}{2\sum_{i\in[N]}(c_i+e_i)^2}\right).$$
    }
    \begin{proof}
    This is Theorem 1.9 of~\cite{Stronger Azuma}.
    \end{proof}
    
    We can now consider the number of isolated $k$-vertex trees in $G[V_p]$:
    \paragraph{Lemma 5.4} \textit{Assume that $\alpha>\frac{2}{\ln(n/d)}$. Let $k\le\frac{1}{\alpha}$. Then, with probability at least $1-\left(\frac{n}{d}\right)^{-\ln\ln(n/d)+8}$, the number of isolated $k$-vertex trees in $G[V_p]$ is at least
    $$(1-4\alpha)\frac{k^{k-2}(1+\epsilon)^k\exp\left(-(1+\epsilon)k\right)}{k!}\cdot\frac{n}{d}.$$
    }
    \begin{proof}
    Fix $M$ satisfying $np-(np)^{2/3}\le M\le np+(np)^{2/3}$ and form $V_M$ by uniformly choosing $M$ vertices from $V$. We will work in $G[V_M]$, use Lemma 5.3, and then convert our results to $G[V_p]$. 
    
    Let $k\le\frac{1}{\alpha}$, and let $T_{k,M}$ be the number of isolated $k$-vertex trees in $G[V_M]$. Let $t_k$ be the number of acyclic connected $k$-sets in $G$. We denote by $(a)_{b}$ the falling factorial $a(a-1)\cdots(a-b+1)$. In order for a $k$-vertex tree in $G$ to become an isolated $k$-vertex tree in $G[V_M]$, we need to include all of its $k$ vertices in $M$, and none of its neighbours, at most $kd$, in $M$. The probability for that is at least
    \begin{align*}
        \frac{{{n-k-kd}\choose {M-k}}}{{n\choose M}}&=\frac{M!}{(M-k)!}\cdot\frac{(n-M)!}{n!}\cdot\frac{(n-k-kd)!}{(n-M-kd)!}\\
        &=\frac{(M)_k}{(n)_M}\cdot\frac{(n-kd)!}{(n-kd-M)!(n-kd)_k}\\
        &\ge\frac{(M)_{k}}{(n)_{k}}\cdot\frac{(n-kd)_M}{(n)_M}\\
        &\ge\left(\frac{np-(np)^{2/3}-k}{n}\right)^k\cdot\left(\frac{n-kd-np-(np)^{2/3}}{n-np-(np)^{2/3}}\right)^{np+(np)^{2/3}},
    \end{align*}
    where we used that $M<n-kd-np-(np)^{2/3}<n$, for our values of $M$ and $k$, in the last inequality. We will treat each multiplicative term separately:
    \begin{align*}
        \left(\frac{np-(np)^{2/3}-k}{n}\right)^k&\ge\left(\frac{1+\epsilon}{d}-\frac{2}{d^{2/3}n^{1/3}}-\frac{k}{n}\right)^k\\
        &\ge\left(\frac{1+\epsilon}{d}\right)^k\left(1-2\left(\frac{d}{n}\right)^{1/3}-\frac{kd}{n}\right)^k\\
        &\ge \left(\frac{1+\epsilon}{d}\right)^k\left(1-2k\left(\left(\frac{d}{n}\right)^{1/3}+\frac{kd}{n}\right)\right)\\
        &\ge\left(\frac{1+\epsilon}{d}\right)^k(1-\alpha),
    \end{align*}
    where we used Bernoulli's Inequality in the third inequality; as for the last inequality, recall that $\alpha \ge \frac{2}{\ln\left(\frac{n}{d}\right)}$ and $k\le \frac{1}{\alpha}$, and hence for $\alpha \ge 2k\left(\left(\frac{d}{n}\right)^{1/3}+\frac{kd}{n}\right)$, it suffices to show that $\alpha^2\ge 3\left(\frac{d}{n}\right)^{1/3}$, and indeed $\frac{4}{\ln^2\left(\frac{n}{d}\right)}\ge \frac{3}{(n/d)^{1/3}}$. Furthermore,
    \begin{align*}
        \left(\frac{n-kd-np-(np)^{2/3}}{n-np-(np)^{2/3}}\right)^{np+(np)^{2/3}}&= \left(1-\frac{kd}{n-np-(np)^{2/3}}\right)^{\frac{(1+\epsilon)n}{d}+\left(\frac{(1+\epsilon)n}{d}\right)^{2/3}}\\
        &\ge\exp\left(-(1+\epsilon)k-\frac{2kd^{1/3}}{n^{1/3}}\right)\\
        &\ge(1-\alpha)\exp\left(-(1+\epsilon)k\right),
    \end{align*}
    where we used that $1-c\ge\exp(-c-c^2)$ for small enough $c$, together with our choice of $\alpha$ and $k$. Together with Lemma 5.2, we obtain:
    \begin{align*}
        \textbf{E}T_{k,m}&\ge(1-2\alpha)t_k\left(\frac{1+\epsilon}{d}\right)^k\exp\left(-(1+\epsilon)k\right)\\
        &\ge(1-3\alpha)\frac{(1+\epsilon)n}{d}\cdot\frac{k^{k-2}(1+\epsilon)^{k-1}\exp\left(-(1+\epsilon)k\right)}{k!}.
    \end{align*}
    We now use the notation of Lemma 5.3. Let $f=T_{k,M}$. Let $\Gamma$ be the event that $\Delta\left(G[V_M]\right)\le\ln\left(\frac{n}{d}\right)$. In order to have a vertex with degree higher than $\ln(n/d)$ in $V_M$, we need to choose the vertex itself and at least $\ln(n/d)$ of its neighbours and include them in $V_M$. Hence:
    \begin{align*}
        P[X\notin\Gamma]&\le n\frac{{d\choose \ln(n/d)}{{n-\ln(n/d)-1}\choose{M-\ln(n/d)}-1}}{{n\choose M}}\\
        &\le n\left(\frac{ed}{\ln(n/d)}\right)^{\ln(n/d)}\frac{(M)_{\ln(n/d)+1}}{(n)_{\ln(n/d)+1}}\\
        &\le \frac{2n}{d}\left(\frac{ed}{\ln(n/d)}\right)^{\ln(n/d)}\frac{(M)_{\ln(n/d)}}{(n)_{\ln(n/d)}}\\
        &\le \frac{2n}{d}\left(\frac{ed}{\ln(n/d)}\cdot\frac{np+(np)^{2/3}}{n}\right)^{\ln(n/d)}\\
        &\le\left(\frac{n}{d}\right)^{-\ln\ln(n/d)+3},
    \end{align*}
    where the third inequality follows from our assumption on $M$. Furthermore, observe that $|e(V_M, V-V_M)|\le dM$. Since there are $n-M$ vertices in $V-V_M$, by the pigeonhole principle there must be a vertex $v_0\in V-V_M$ such that $d(v_0,V_M)\le 3$. We define $\rho_i$ be the bijection such that $\rho_i(x)$ satisfies that its $i$-th entry is $v_0$ and all the other entries remain unchanged. Since $V_M$ is chosen uniformly among all sets of size $M$ in $V$, the condition of Lemma 5.3 is satisfied: $$P[X=x|X\in\Sigma_a]=\frac{P[X=x]}{P[X\in\Sigma_a]}=\frac{P[X=x]}{P[X\in\Sigma_{v_0}]}=P[X=x|X\in\Sigma_{v_0}].$$
    Considering $\Big|f(x)-f\left(\rho_i(x)\right)\Big|$, we need only to consider the possible change in the value of $T_{k,M}$ when we change one vertex. Clearly, it cannot change by more than $M$. Furthermore, If $x$ is in $\Gamma$, $\Delta(G[V_M])\le\ln(n/d)$, and since $\rho_i(x)$ chooses $v_0$ which has less than $\ln(n/d)$ neighbours in $V_M$, this change in one vertex cannot change the value of $T_{k,M}$ by more than $\ln(n/d)$. As such, the conditions of Lemma 5.3 hold with $c_i=\ln(n/d)$ and $d_i=M$. Choosing $\gamma_i=\frac{1}{M}$, we have that $e_i\le 1$ and
    \begin{align*}
        P[\mathcal{B}]&\le\sum_{i=1}^MM\cdot P[X\notin\Gamma]\\
        &\le M^2\left(\frac{n}{d}\right)^{-\ln\ln(n/d)+3}\\
        &\le \left(\frac{n}{d}\right)^{-\ln\ln(n/d)+6}.
    \end{align*}
    We thus obtain by Lemma 5.3:
    \begin{align*}
        P&\left[T_{k,M}\le \frac{(1-4\alpha)nk^{k-2}\left((1+\epsilon)\exp\left(-(1+\epsilon)\right)\right)^{k}}{k!d}\right]\\
        &\le P[\neg\mathcal{B}]+ \exp\left(-\frac{(np)^{4/3}}{4\sum_{i=1}^M\left(\ln(n/d)+1\right)^2}\right)\\
        &\le \left(\frac{n}{d}\right)^{-\ln\ln(n/d)+6}+\exp\left(-\left(\frac{n}{d}\right)^{1/4}\right)\\
        &\le\left(\frac{n}{d}\right)^{-\ln\ln(n/d)+7},
    \end{align*}
    where we note that $(np)^{2/3}<\alpha\frac{nk^{k-2}(1+\epsilon)^{k}\exp\left(-(1+\epsilon)k\right)}{k!d}$ for our choice of $\alpha$ and $k$.
    
    We now convert our results to $G[V_p]$. Let $T_{k,p}$ be the number of isolated $k$-vertex trees in $G[V_p]$. Note that by Chernoff,
    $$P\left[\Big||V_p|-np\big|\ge(np)^{2/3}\right]\le \exp\left(-\left(\frac{n}{d}\right)^{1/4}\right).$$
    Let $\mathcal{P}$ be an arbitrary graph property. By the Law of Total Probability,
    \begin{align*}
        P\left[G[V_p]\in\mathcal{P}\right]&=\sum_{m=0}^nP[|V_p|=M]\cdot P\left[G[V_M]\in\mathcal{P}\right]\\
        &\le \exp\left(-\left(\frac{n}{d}\right)^{1/4}\right)+\sum_{M=np-(np)^{2/3}}^{np+(np)^{2/3}}P[|V_p|=M]\cdot P\left[G[V_M]\in\mathcal{P}\right].
    \end{align*}
    We may thus conclude:
    \begin{align*}
        P\left[T_{k,p}\le (1-4\alpha)\frac{nk^{k-2}(1+\epsilon)^{k}\exp\left(-(1+\epsilon)k\right)}{k!d}\right]\le \left(\frac{n}{d}\right)^{-\ln\ln(n/d)+8}.
    \end{align*}
    \end{proof}

    We are now ready to prove Theorem 3.
    \begin{proof}[\textbf{\textit{Proof of Theorem 3.}}]
    Let $k\le\frac{1}{\alpha}$. Let $\beta_1$, $\beta_2$ and $\beta_3$ be some positive constants (possibly depending on $\epsilon$) to be determined later. Let $T_{k,p}$ denote the number of isolated $k$-vertex trees in $G[V_p]$. By Lemma 5.4, with probability at least $1-\left(\frac{n}{d}\right)^{-\ln\ln(n/d)+8}$, we have that
    $$T_{k,p}\ge(1-4\alpha)\frac{k^{k-2}(1+\epsilon)^{k}\exp\left(-(1+\epsilon)k\right)}{k!}\cdot\frac{n}{d}.$$ 
    Thus, with probability at least $1-\left(\frac{n}{d}\right)^{-\beta_1\ln\ln(n/d)}$, the number of vertices in isolated trees of order at most $\frac{1}{\alpha}$ is at least:
    \begin{align*}
        \sum_{k=1}^{\frac{1}{\alpha}}kT_{k,p}&\ge (1-4\alpha)\frac{(1+\epsilon)n}{d}\sum_{k=1}^{\frac{1}{\alpha}}\frac{k^{k-1}}{k!}(1+\epsilon)^{k-1}\exp\left(-(1+\epsilon)k\right).
    \end{align*}
    Now, we want to bound the sum $$\sum_{j=\frac{1}{\alpha}}^{\infty}\frac{j^{j-1}}{j!}(1+\epsilon)^{j}\exp\left(-(1+\epsilon)j\right).$$
    For that, define
    $$f(j)=\frac{j^{j-1}}{j!}(1+\epsilon)^j\exp\left(-(1+\epsilon)j\right),$$
    with $j\in\mathbb{N}, j\ge\frac{1}{\alpha}$. Observe that:
    \begin{align*}
        \frac{f(j+1)}{f(j)}&=(1+\epsilon)\exp(-1-\epsilon)\frac{(j+1)^j}{j^{j-1}}\cdot\frac{1}{j+1}\\
        &\le (1+\epsilon)\left(1-\epsilon+\frac{\epsilon^2}{2}\right)\exp(-1)\left(\frac{j+1}{j}\right)^{j-1}\\
        &\le 1-\frac{\epsilon^2}{3},
    \end{align*}
    where we used inequality $(7)$ in the first inequality. 
    Thus, for all $j\ge \frac{1}{\alpha}$, we have that 
    $$f(j)\le f\left(\frac{1}{\alpha}\right)\left(1-\frac{\epsilon^2}{3}\right)^{j-\frac{1}{\alpha}}.$$
    As such,
    \begin{align*}
        \sum_{j=\frac{1}{\alpha}}^{\infty}\frac{j^{j-1}}{j!}(1+\epsilon)^{j}\exp\left(-(1+\epsilon)j\right)&=\sum_{j=\frac{1}{\alpha}}^{\infty}f(j)\\
        &\le \frac{3f\left(\frac{1}{\alpha}\right)}{\epsilon^2}.
    \end{align*}
    Using the inequality $n!>n^{n+1/2}\exp(-n)$, we obtain:
    \begin{align*}
        f\left(\frac{1}{\alpha}\right)&\le \alpha^{3/2}\exp\left(\frac{1}{\alpha}\right)\exp\left(-(1+\epsilon)\frac{1}{\alpha}\right)\\
        &\le \alpha^{3/2}.
    \end{align*}
    We may thus conclude,
    \begin{align*}
         \sum_{j=\frac{1}{\alpha}}^{\infty}\frac{j^{j-1}}{j!}(1+\epsilon)^{j}\exp\left(-(1+\epsilon)j\right)&\le \frac{3\alpha^{3/2}}{\epsilon^2}\\
         &\le 3\alpha,
    \end{align*}
    since we assume $\alpha<\epsilon^4$. We may thus conclude:
    \begin{align*}
       \sum_{k=1}^{\frac{1}{\alpha}}kT_{k,p}&\ge (1-4\alpha)\frac{(1+\epsilon)n}{d}\sum_{k=1}^{\frac{1}{\alpha}}\frac{k^{k-1}}{k!}(1+\epsilon)^{k-1}\exp\left(-(1+\epsilon)k\right)\\
       &\ge (1-4\alpha)\frac{(1+\epsilon)n}{d}\sum_{k=1}^{\infty}\frac{k^{k-1}}{k!}(1+\epsilon)^{k-1}\exp\left(-(1+\epsilon)k\right)-\frac{3\alpha n}{d}.
    \end{align*}
    It is known (see, for example,~\cite{Erdos1}) that:
    $$\sum_{k=1}^\infty \frac{k^{k-1}}{k!}(1+\epsilon)^{k-1}\exp\left(-(1+\epsilon)k\right)=\frac{y}{1+\epsilon},$$
    where $y$ is as defined in $(3)$. Using this together with $(4)$, we obtain that with probability at least $1-\left(\frac{n}{d}\right)^{-\beta_1\ln\ln(n/d)}$, we have at least
    $$\frac{yn}{d}-\frac{7\alpha n}{d}=\frac{(1+\epsilon-x)n}{d}-\frac{7\alpha n}{d}$$
    vertices in isolated trees of order at most $\frac{1}{\alpha}$ in $G[V_p]$. 
    By Theorem 2, we have with probability at least $1-\exp\left(-\beta_2\frac{\alpha^2n}{d}\right)$ at least $\frac{xn}{d}-\frac{7\alpha n}{d}$ vertices in a large component. By Property 3 of Lemma 2.1 (using $\alpha$ instead of $c$), we have that with probability at least $1-\exp\left(-\beta_3\frac{\alpha^2 n}{d}\right)$ the number of vertices in $G[V_p]$ is at most
    $$\sum_{i=1}^nX_i\le\frac{(1+\epsilon)n}{d}+\frac{\epsilon^2\alpha n}{d}.$$
    We thus conclude that with probability at least $1-\left(\frac{n}{d}\right)^{-\beta\ln\ln(n/d)}$, there is only one large component, and there are at most $$\frac{(1+\epsilon)n}{d}+\frac{\epsilon^2\alpha n}{d}-\left(\frac{(1+\epsilon-x)n}{d}-\frac{7\alpha n}{d}\right)-\left(\frac{xn}{d}-\frac{7\alpha n}{d}\right)\le\frac{15\alpha n}{d}$$ vertices which are neither in the large component nor in isolated trees of order at most $\frac{1}{\alpha}$.
    \end{proof}
    Theorems 2 and 3 together provide us with a description of the components of $G[V_p]$ in the supercritical phase (with superpolynomially high probability): we conclude that there is a unique largest component, and its size is concentrated around $\frac{xn}{d}$, where $x$ is as defined by $(2)$. Furthermore, the second largest component is of size $O(\alpha n/d)$, and the volume of vertices in isolated trees of small size is concentrated around $\frac{yn}{d}$, where $y$ is as defined by $(3)$. 
    
\section{Further Properties of the Giant Component}
    First, we prove a tight concentration result for the number of edges in $G[V_p]$:
    
    \paragraph{Lemma 6.1} \textit{Let $Z_p$ be the random variable counting the number of edges in $G[V_p]$. Then, for some constant $\beta>0$, with probability at least $1-\left(\frac{n}{d}\right)^{-\beta\ln\ln(n/d)}$,
    $$\Bigg|Z_p-\frac{(1+\epsilon)^2n}{2d}\Bigg|\le 2(np)^{2/3}.$$
    }
    \begin{proof} 
    Fix $M$ satisfying $np-(np)^{2/3}\le M\le np+(np)^{2/3}$ and form $V_M$ by uniformly choosing $M$ vertices from $V$. Let $Z_M$ be the random variable representing the number of edges in $G[V_M]$. As in Lemma 5.4, we will work in $G[V_M]$, use Lemma 5.3, and then convert our results to $G[V_p]$. 
    
    For an edge $uv$ to be included in $G[V_M]$, we need to include both $u$ and $v$ in $V_M$, which happens with probability:
    $$\frac{{{n-2}\choose {M-2}}}{{n\choose M}}=\frac{M(M-1)}{n(n-1)}.$$
    Since $G$ is a $d$-regular graph, it has $\frac{nd}{2}$ edges. As such,
    \begin{align*}
    \textbf{E}Z_M=\frac{M(M-1)d}{2(n-1)}.
    \end{align*}
    As in Lemma 5.4, we now use the notation of Lemma 5.3. Let $f=Z_M$. Let $\Gamma$ be the event that $\Delta\left(G[V_M]\right)\le\ln\left(\frac{n}{d}\right)$.Furthermore, observe that $|e(V_M, V-V_M)|\le dM$. Since there are $n-M$ vertices in $V-V_M$, by the pigeonhole principle there must be a vertex $v_0\in V-V_M$ such that $d(v_0,V_M)\le 3$. We define $\rho_i$ be the bijection such that $\rho_i(x)$ satisfies that its $i$-th entry is $v_0$ and all the other entries remain unchanged. Since $V_M$ is chosen uniformly among all sets of size $M$ in $V$, the condition of Lemma 5.3 is satisfied: $$P[X=x|X\in\Sigma_a]=\frac{P[X=x]}{P[X\in\Sigma_a]}=\frac{P[X=x]}{P[X\in\Sigma_{v_0}]}=P[X=x|X\in\Sigma_{v_0}].$$
    Considering $\Big|f(x)-f\left(\rho_i(x)\right)\Big|$, we need only to consider the possible change in the value of $f$ when we change one vertex. As in Lemma 5.4, it cannot change by more than $M$, and when $x$ is in $\Gamma$, it cannot change by more than $\ln(n/d)$. As such, the conditions of Lemma 5.3 hold with $c_i=\ln(n/d)$ and $d_i=M$. Thus, as in Lemma 5.4, we choose $\gamma_i=\frac{1}{M}$, and we obtain that $P[\mathcal{B}]\le\left(\frac{n}{d}\right)^{-\ln\ln(n/d)+6}$, and
    $$P\left[Z_M\le\textbf{E}Z_M-(np)^{2/3}\right]\le\left(\frac{n}{d}\right)^{-\ln\ln(n/d)+7}.$$
    Similarly, defining $g=-f$, we obtain that:
    $$P\left[Z_M\ge\textbf{E}Z_M+(np)^{2/3}\right]\le\left(\frac{n}{d}\right)^{-\ln\ln(n/d)+7}.$$
    As in Lemma 5.4, we can convert these results from $G[V_M]$, where $np-(np)^{2/3}\le M\le np+(np)^{2/3}$, and conclude that:
    $$P\left[\Bigg|Z_p-\frac{(1+\epsilon)^2n}{2d}\Bigg|\le (np)^{2/3}\right]\le\left(\frac{n}{d}\right)^{-\beta\ln\ln(n/d)}.$$
    To complete the proof note that $\Bigg|\frac{M(M-1)d}{2(n-1)}-\frac{(1+\epsilon)^2n}{2d}\Bigg|\le (np)^{2/3}$, for our values of $M$.
    \end{proof}
    
    We now proceed to bound the number of edges in trees of order at most $\frac{1}{\sqrt{\alpha}}$. 
    \paragraph{Lemma 6.2} \textit{Assume that $\frac{2}{\ln(n/d)}<\alpha<\epsilon^8$.Let $Z_p$ be the random variable counting the number of edges in isolated $k$-vertex trees in $G[V_p]$, where $1\le k\le \frac{1}{\sqrt{\alpha}}$. Then, there exists some positive constant $\beta$ such that with probability at least $1-\left(\frac{n}{d}\right)^{-\beta\ln\ln(n/d)}$, 
    $$\Bigg|Z_p-\frac{(1+\epsilon-x)^2n}{2d}\Bigg|\le \frac{16\sqrt{\alpha}n}{d},$$
    where $x$ is as defined in $(2)$.
    }
    \begin{proof}
    Let $\beta_1$ through $\beta_5$ be some positive constants (possibly depending on $\epsilon$) to be determined later. Denote by $T_{k,p}$ the number of isolated $k$-vertex trees in $G[V_p]$. By Lemma 5.4, with probability at least $1-\left(\frac{n}{d}\right)^{-\ln\ln(n/d)+8}$,
    $$T_{k,p}\ge(1-4\alpha)\frac{n}{d}\cdot\frac{k^{k-2}(1+\epsilon)^k\exp\left(-(1+\epsilon)k\right)}{k!}.$$ 
    As such, with probability at least $1-\left(\frac{n}{d}\right)^{-\beta_1\ln\ln(n/d)}$,
    \begin{align*}
        Z_p&\ge(1-4\alpha)\frac{n}{d}\sum_{k=1}^{\frac{1}{\sqrt{\alpha}}}\frac{(k-1)k^{k-2}(1+\epsilon)^k\exp\left(-(1+\epsilon)k\right)}{k!}.
    \end{align*}
    Now, we want to bound the sum $$\sum_{j>\frac{1}{\sqrt{\alpha}}}\frac{(j-1)j^{j-2}(1+\epsilon)^j\exp\left(-(1+\epsilon)j\right)}{j!}.$$
    For that, similar the proof of Theorem 3, we define
    $$f(j)=\frac{(j-1)j^{j-2}(1+\epsilon)^j\exp\left(-(1+\epsilon)j\right)}{j!},$$
    with $j\in\mathbb{N}, j\ge \frac{1}{\sqrt\alpha}$. Observe that:
    \begin{align*}
        \frac{f(j+1)}{f(j)}&=(1+\epsilon)\exp\left(-1-\epsilon\right)\frac{j}{j-1}\cdot\frac{(j+1)^{j-1}}{j^{j-2}}\cdot\frac{1}{j+1}\\
        &\le (1+\epsilon)\left(1-\epsilon+\frac{\epsilon^2}{2}\right)\exp(-1)\left(\frac{j}{j-1}\right)^{j-1}\\
        &\le 1-\frac{\epsilon^2}{3},
    \end{align*}
    where we used inequality $(7)$ in the first inequality.
    Thus, for all $j\ge\frac{1}{\alpha}$ we have that
    $$f(j)\le f\left(\frac{1}{\alpha}\right)\left(1-\frac{\epsilon^2}{3}\right)^{j-\frac{1}{\alpha}}.$$
    As such,
    \begin{align*}
        \sum_{j=\frac{1}{\sqrt{\alpha}}}^{\infty}\frac{(j-1)j^{j-2}(1+\epsilon)^j\exp\left(-(1+\epsilon)j\right)}{j!}&=\sum_{j=\frac{1}{\sqrt{\alpha}}}^{\infty}f(j)\\
        &\le\frac{3f\left(\frac{1}{\sqrt{\alpha}}\right)}{\epsilon^2}.
    \end{align*}
    Using the inequality $n!<n^{n+1/2}\exp(-n)$, we obtain:
    \begin{align*}
        f\left(\frac{1}{\sqrt{\alpha}}\right)\le \alpha^{3/4}.
    \end{align*}
    We may thus conclude,
    \begin{align*}
        \sum_{j=\frac{1}{\sqrt{\alpha}}}^{\infty}\frac{(j-1)j^{j-2}(1+\epsilon)^j\exp\left(-(1+\epsilon)j\right)}{j!}&\le \frac{3\alpha^{3/4}}{\epsilon^2}\\
        &\le 3\sqrt{\alpha},
    \end{align*}
    since we assume $\alpha<\epsilon^8$. We may now conclude that with probability at least $1-\left(\frac{n}{d}\right)^{-\beta_1\ln\ln(n/d)}$:
    \begin{align*}
         Z_p&\ge(1-4\alpha)\frac{n}{d}\sum_{k=1}^{\infty}\frac{(k-1)k^{k-2}(1+\epsilon)^k\exp\left(-(1+\epsilon)k\right)}{k!}-\frac{3\sqrt{\alpha}n}{d}.
    \end{align*}
    It is known (see, for example,  Theorem 2.14 in~\cite{Intro to Random Graphs} and its proof) that the number of edges in $G\left(n,\frac{1+\epsilon}{n}\right)$ which lie in isolated trees is asymptotically $\frac{y^2n}{2(1+\epsilon)}$, and as such:
    $$\sum_{k=1}^\infty\frac{(k-1)k^{k-2}}{k!}\left((1+\epsilon)\exp\left(-(1+\epsilon)\right)\right)^k=\frac{y^2}{2},$$
    where $y$ is as defined in $(3)$. Using this together with $(4)$, we obtain that with probability at least $1-\left(\frac{n}{d}\right)^{-\beta_1\ln\ln(n/d)}$
    \begin{align*}
        Z_p&\ge (1-4\alpha)\frac{n}{d}\sum_{k=1}^{\infty}\frac{(k-1)k^{k-2}(1+\epsilon)^k\exp\left(-(1+\epsilon)k\right)}{k!}-\frac{3\sqrt{\alpha}n}{d}\\
        &\ge \frac{y^2n}{2d}-\frac{4\sqrt{\alpha}n}{d}\\
        &\ge \frac{(1+\epsilon-x)^2n}{2d}-\frac{4\sqrt{\alpha}n}{d}.
    \end{align*}

    For the other side, recall that by Property 3 of Lemma 2.1, with probability at least $1-\exp\left(-\beta_2\frac{\alpha^2 n}{d}\right)$ we have that
    $$|V_p|<\frac{(1+\epsilon+\alpha)n}{d}.$$
    Furthermore, by Theorem 2, with probability at least $1-\exp\left(-\beta_3\frac{\alpha^2 n}{d}\right)$ there is a component of size at least $\frac{(x-7\alpha)n}{d}$ in $V_p$.
    Now, assume for contradiction that for some $k$, $1\le k \le \frac{1}{\sqrt{\alpha}}$, we have that with probability at least $\left(\frac{n}{d}\right)^{-\beta_4\ln\ln(n/d)}$ that
    $$T_{k,p}\ge\frac{16\alpha n}{kd}+\frac{n}{d}\cdot\frac{k^{k-2}(1+\epsilon)^{k-1}\exp\left(-(1+\epsilon)k\right)}{k!}.$$
    We will show that this contradicts Property 3 of Lemma 2.1, that is, too high probability $V_p$ will be too large. Indeed, together with Lemma 5.4 and the union bound, we can conclude that with probability at least $\frac{2}{3}\left(\frac{n}{d}\right)^{-\beta_4\ln\ln(n/d)}$, the number of vertices in isolated trees of order at most $\frac{1}{\alpha}$ is at least
    \begin{align*}
        &\frac{16\alpha n}{d}+(1-4\alpha)\frac{(1+\epsilon)n}{d}\cdot\sum_{k=1}^{\frac{1}{\alpha}}\frac{k^{k-1}(1+\epsilon)^{k-1}\exp\left(-(1+\epsilon)k\right)}{k!}\\
        &\ge
        \frac{16\alpha n}{d}+(1-7\alpha)\frac{(1+\epsilon-x)n}{d},
    \end{align*}
    where we used the same lower bound for the sum as in the proof of Theorem 3. Altogether, we conclude that with probability at least $\frac{1}{2}\left(\frac{n}{d}\right)^{-\beta_4\ln\ln(n/d)}$,
    \begin{align*}
        |V_p|&\ge\frac{(x-7\alpha)n}{d}+\frac{16\alpha n}{d}+(1-7\alpha)\frac{(1+\epsilon-x)n}{d}\\
        &>\frac{(1+\epsilon+\alpha)n}{d},
    \end{align*}
    which is a contradiction to Property 3 of Lemma 2.1 if $$\frac{1}{2}\left(\frac{n}{d}\right)^{-\beta_4\ln\ln(n/d)}>\exp\left(-\beta_3\frac{\alpha^2n}{d}\right),$$ 
    that is with the right choice of $\beta_4$. Thus, with probability at least $1-\left(\frac{n}{d}\right)^{-\beta_5\ln\ln(n/d)}$, we have that for all $1\le k \le \frac{1}{\sqrt{\alpha}}$:
    \begin{align*}
        T_{k,p}\le \frac{16\alpha n}{kd}+\frac{n}{d}\cdot\frac{k^{k-2}(1+\epsilon)^{k-1}\exp\left(-(1+\epsilon)k\right)}{k!}.
    \end{align*}
    As such, with probability at least $1-\left(\frac{n}{d}\right)^{-\beta_7\ln\ln(n/d)}$,
    \begin{align*}
        Z_p&=\sum_{k=1}^{\frac{1}{\sqrt{\alpha}}}(k-1)T_{k,p}\\
        &\le\sum_{k=1}^{\frac{1}{\sqrt{\alpha}}}\left(\frac{16\alpha n}{d}+\frac{n}{d}\cdot\frac{(k-1)k^{k-2}(1+\epsilon)^{k-1}\exp\left(-(1+\epsilon)k\right)}{k!}\right) \\
        &\le \frac{16\sqrt{\alpha}n}{d}+\frac{n}{d}\cdot\sum_{k=1}^\infty\frac{(k-1)k^{k-2}}{k!}\left((1+\epsilon)\exp\left(-(1+\epsilon)\right)\right)^k\\
        &\le \frac{(1+\epsilon-x)^2n}{2d}+\frac{16\sqrt{\alpha}n}{d},
    \end{align*}
    completing the proof. 
    \end{proof}
    
    In order to obtain our bound on the number of edges in the giant component, we need to bound the number of edges in components which are neither the giant component, nor isolated trees of order at most $\frac{1}{\sqrt{\alpha}}$. We first bound the number of vertices in such components:
    \paragraph{Lemma 6.3} \emph{Assume that $\frac{2}{\ln(n/d)}<\alpha<\epsilon^8$. Then, there exists a positive constant $\beta$ such that with probability at least $1-\left(\frac{n}{d}\right)^{-\beta\ln\ln(n/d)}$, the number of vertices neither in the giant component and nor in isolated trees of order at most $\frac{1}{\sqrt{\alpha}}$ is at most $\frac{5\sqrt{\alpha}n}{d}$.
    }
    \begin{proof}
    Let $\beta_1$, $\beta_2$, $\beta_3$ be some positive constants (possibly depending on $\epsilon$) to be determined later. By Property 3 of Lemma 2.1, with probability at least $1-\exp\left(-\beta_1\frac{\alpha^2 n}{d}\right)$ we have that
    $$|V_p|<\frac{(1+\epsilon+\alpha)n}{d}.$$
    Furthermore, by Theorem 2, with probability at least $1-\exp\left(-\beta_2\frac{\alpha^2 n}{d}\right)$ there is a component of order at least $\frac{(x-7\alpha)n}{d}$ in $V_p$.
    By Lemma 5.4 together with the union bound, we have that with probability at least $1-\left(\frac{n}{d}\right)^{-\beta_3\ln\ln(n/d)}$, the number of vertices in trees of order at most $\frac{1}{\sqrt{\alpha}}$ is at least:
    \begin{align*}
        \sum_{k=1}^{\frac{1}{\sqrt{\alpha}}}kT_{k,p}&\ge (1-4\alpha)\frac{(1+\epsilon)n}{d}\sum_{k=1}^{\frac{1}{\sqrt{\alpha}}}\frac{k^{k-1}(1+\epsilon)^k\exp\left(-(1+\epsilon)k\right)}{k!}\\
        &\ge (1-4\alpha)\frac{(1+\epsilon)n}{d}\sum_{k=1}^{\infty}\frac{k^{k-1}(1+\epsilon)^k\exp\left(-(1+\epsilon)k\right)}{k!}-\frac{3\sqrt{\alpha}n}{d}\\
        &\ge \frac{(1+\epsilon-x)n}{d}-\frac{4\sqrt{\alpha}n}{d},
    \end{align*}
    where we used our bound on the tail of the sum from the proof of Lemma 6.2, and our bound on the sum itself from the proof of Theorem 3. Altogether, we have with probability at least $1-\left(\frac{n}{d}\right)^{-\beta\ln\ln(n/d)}$, at most
    \begin{align*}
        \frac{(1+\epsilon+\alpha)n}{d}-\frac{(x-7\alpha)n}{d}-\frac{(1+\epsilon-x)n}{d}+\frac{4\sqrt{\alpha}n}{d}\le \frac{5\sqrt{\alpha}n}{d}
    \end{align*}
    vertices which are neither in the giant component nor in isolated trees of order at most $\frac{1}{\sqrt{\alpha}}$.
    \end{proof}
    
    We may now bound the number of edges in such components:
    \paragraph{Lemma 6.4}\emph{Assume that $\frac{2}{\ln(n/d)}<\alpha<\epsilon^8$. Then, there exists a positive constant $\beta$ such that with probability at least $1-\left(\frac{n}{d}\right)^{-\beta\ln\ln(n/d)}$, the number of edges in $G[V_p]$ which are in components that are neither the giant component nor isolated trees of order at most $\frac{1}{\sqrt{\alpha}}$ is at most $\frac{7\alpha^{1/4}n}{d}$.}
    \begin{proof}
    Let $\beta_1,\beta_2$ be positive constants (possibly depending on $\epsilon$) to be determined later. Let $k_0=\frac{1}{\alpha^{1/4}}$. Fix $M$ satisfying $np-(np)^{2/3}\le M\le np+(np)^{2/3}$ and form $V_M$ by uniformly choosing $M$ vertices from $V$. As in Lemma 5.4 and 6.1, we will work in $G[V_M]$, use Lemma 5.3, and then convert our results to $G[V_p]$.
    
    Let $Z_M$ be the number of edges in components that are neither the giant component nor isolated trees of order at most $\frac{1}{\sqrt{\alpha}}$ in $G[V_M]$. Denote the set of vertices which are neither in the giant component nor in isolated trees of order at most $\frac{1}{\sqrt{\alpha}}$ by $S$. Let $$D_{k_0,M}=\{v\in V_M: d_{V_M}(v)>k_0\}.$$ We define the random variable 
    $$Y=\sum_{v\in D_{k_0,M}}d_{V_M}(v).$$
    By Lemma 6.3, with probability at least $1-\left(\frac{n}{d}\right)^{-\beta_1\ln\ln(n/d)}$,
    we have at most $\frac{5\sqrt{\alpha}n}{d}$ vertices in $S$. We thus obtain:
    $$\textbf{E}Z_M\le \frac{5\sqrt{\alpha}n}{d}\cdot k_0+\left(\frac{n}{d}\right)^{-\beta_1\ln\ln(n/d)}\cdot M+\textbf{E}Y.$$
    Observe that:
    \begin{align*}
        \textbf{E}Y\le n\cdot\sum_{i=k_0}^{M}i\cdot P[d_{V_M}(v)=i],
    \end{align*}
    since we have $n$ vertices to consider, and then consider the expected degree of each of these vertices in $V_M$ (which is at least $k_0$ and at most $M$). Furthermore, 
    \begin{align*}
        \sum_{i=k_0}^M i\cdot P[d_{V_M}(v)=i]&=\sum_{i=k_0}^M i\cdot \left(P[d_{V_M}(v)\ge i]-P[d_{V_M}(v)\ge i+1]\right)\\
        &\le \sum_{i=k_0}^M P[d_{V_M}(v)\ge i]. 
    \end{align*}
    To calculate the probability a vertex has degree at least $i$ in $V_M$, we need to choose the vertex itself and include it in $V_M$, and choose at least $i$ of its $d$ neighbours and include them in $V_M$:
    \begin{align*}
      \frac{{{d}\choose{i}}{{n-i-1}\choose{M-i-1}}}{{n\choose M}}&\le\left(\frac{ed}{i}\right)^{i}\frac{(M)_{i+1}}{(n)_{i+1}} \\
      &\le \frac{2}{d}\left(\frac{ed}{i}\right)^{i}\frac{(M)_{i}}{(n)_{i}}\\
      &\le\frac{2}{d}\left(\frac{ed}{i}\cdot\frac{np+(np)^{2/3}}{n}\right)^{i}\\
      &\le\frac{2}{d}\left(\frac{2e}{i}\right)^{i}.
    \end{align*}
    Recall that $k_0=\frac{1}{\alpha^{1/4}}$. We thus obtain:
    \begin{align*}
        \textbf{E}Y&\le n\cdot\sum_{i=k_0}^{M}i\cdot P[d_{V_M}(v)=i]\\
        &\le n\cdot\sum_{i=k_0}^M P[d_{V_M}(v)\ge i] \\
        &\le n\cdot\sum_{i=k_0}^M \frac{2}{d}\left(\frac{2e}{i}\right)^i\\
        &\le \frac{2n}{d}\sum_{i=k_0}^M \left(\frac{2e}{k_0}\right)^i \\
        &\le \frac{2n}{d}(3e\alpha^{1/4})^{\frac{1}{\alpha^{1/4}}},
    \end{align*}
    where we used the fact that 
    $$\sum_{i=a}^{b}\left(\frac{2e}{a}\right)^i=\frac{\left(\frac{1}{a}\right)^{a-1}(2e)^{a}-(2e)^{b+1}\left(\frac{1}{a}\right)^b}{a-2e}.$$
    We can now conclude:
    \begin{align*}
        \textbf{E}Z_M&\le \frac{5\sqrt{\alpha}n}{d}k_0+\left(\frac{n}{d}\right)^{-\beta_1\ln\ln(n/d)}+\textbf{E}Y\\
        &\le \frac{5\alpha^{1/4}n}{d}+\frac{3n}{d}(3e\alpha^{1/4})^{\frac{1}{\alpha^{1/4}}}\\
        &\le \frac{6\alpha^{1/4}n}{d}.
    \end{align*}

    Now, as in Lemma 5.4 and 6.1, we use the notation of Lemma 5.3. Let $f=Z_M$. Let $\Gamma$ be the event that $\Delta\left(G[V_M]\right)\le\ln\left(\frac{n}{d}\right)$. 
    Furthermore, observe that $|e(V_M, V-V_M)|\le dM$. Since there are $n-M$ vertices in $V-V_M$, by the pigeonhole principle there must be a vertex $v_0\in V-V_M$ such that $d(v_0,V_M)\le 3$. We define $\rho_i$ be the bijection such that $\rho_i(x)$ satisfies that its $i$-th entry is $v_0$ and all the other entries remain unchanged. Since $V_M$ is chosen uniformly among all sets of size $M$ in $V$, the condition of Lemma 5.3 is satisfied: $$P[X=x|X\in\Sigma_a]=\frac{P[X=x]}{P[X\in\Sigma_a]}=\frac{P[X=x]}{P[X\in\Sigma_{v_0}]}=P[X=x|X\in\Sigma_{v_0}].$$
    Considering $\Big|f(x)-f\left(\rho_i(x)\right)\Big|$, changing one vertex cannot change $f$ by more than $M$, and if $x$ is in $\Gamma$, it cannot change $f$ by more than $\ln(n/d)$. As such, the conditions of Lemma 5.3 hold with $c_i=\ln(n/d)$ and $d_i=M$. As in Lemma 5.4, we can choose $\gamma_i=\frac{1}{M}$, to conclude with probability at least $1-\left(\frac{n}{d}\right)^{-\beta_2\ln\ln(n/d)}$,
    $Z_M\le \frac{7\alpha^{1/4}n}{d}$. 
    
    As in Lemma 5.4, since this result holds for all $M$ in the range $np-(np)^{2/3}\le M\le np+(np)^{2/3}$, we can conclude that with probability at least $1-\left(\frac{n}{d}\right)^{-\beta\ln\ln(n/d)}$ the number of edges in $G[V_p]$ in that are neither in the giant component nor on isolated trees of order at most $\frac{1}{\sqrt{\alpha}}$ is at most $\frac{7\alpha^{1/4}n}{d}$.
    \end{proof}

    We are now ready to prove Theorem 4:
    \begin{proof}[\textbf{\textit{Proof of Theorem 4}}]
    We begin with the lower bound. Let $\beta_1$ through $\beta_5$ be positive constants (possibly depending on $\epsilon$) to be determined later. By Lemma 6.1, the number of edges in $G[V_p]$ is, with probability at least $1-\left(\frac{n}{d}\right)^{-\beta_1\ln\ln(n/d)}$, at least 
    $$\frac{(1+\epsilon)^2n}{2d}-\frac{\alpha n}{d}.$$
    By Lemma 6.2, there are, with probability at least $1-\left(\frac{n}{d}\right)^{-\beta_2\ln\ln(n/d)}$, at most 
    $$\frac{(1+\epsilon-x)^2n}{2d}+\frac{16\sqrt{\alpha}n}{d},$$
    edges in isolated trees of order at most $\frac{1}{\sqrt{\alpha}}$ in $G[V_p]$. By Lemma 6.4, besides these edges, the other edges not in the giant component are, with probability at least $1-\left(\frac{n}{d}\right)^{-\beta_3\ln\ln(n/d)}$, at most:
    $$\frac{7\alpha^{1/4}n}{d}.$$
    We thus obtain, that with probability at least $1-\left(\frac{n}{d}\right)^{-\beta_4\ln\ln(n/d)}$, the number of edges in the giant component is at least:
    \begin{align*}
    \frac{(1+\epsilon)^2n}{2d}-\frac{(1+\epsilon-x)^2n}{2d}-\frac{8\alpha^{1/4}}{d}.
    \end{align*}
    
    For the upper bound, we have by Lemma 6.1 that the number of edges in $G[V_p]$, with probability at least $1-\left(\frac{n}{d}\right)^{-\beta_1\ln\ln(n/d)}$, is at most
    $$\frac{(1+\epsilon)^2n}{2d}+\frac{\alpha n}{d}.$$
    By Lemma 6.2, we have that with probability at least $1-\left(\frac{n}{d}\right)^{-\beta_2\ln\ln(n/d)}$, the number of edges in isolated trees of order at most $\frac{1}{\sqrt{\alpha}}$ is at least
    $$\frac{(1+\epsilon-x)^2n}{2d}-\frac{16\sqrt{\alpha}n}{d}.$$
    Thus, with probability at least $1-\left(\frac{n}{d}\right)^{-\beta_5\ln\ln(n/d)}$, the number of edges in the giant component is at most:
    \begin{align*}
    \frac{(1+\epsilon)^2n}{2d}-\frac{(1+\epsilon-x)^2n}{2d}+\frac{20\sqrt{\alpha}n}{d}.
    \end{align*}
    \end{proof}
    
    In~\cite{site-percolation}, Krivelevich proved that the $G[V_p]$ contains \textbf{whp} a path of length $\frac{\epsilon^2n}{5d}$. As a matter of fact, it is easy to see that this result actually holds with probability $1-\exp\left(-\beta\alpha^2 n/d\right)$ for some positive constant $\beta(\epsilon)$. In order to show the existence of a long cycle, we will use Janson's inequality, specifically the following statement from~\cite{Chernoff from K and S}:
    \paragraph{Lemma 6.5} \textit{Let $G=(V,E)$ be a graph with $E=\{e_1, \cdots, e_m\}$. Let $A_i$ be the event that $e_i$ is in $G[V_p]$. Let $X=X_1+\cdots+X_m$, where $X_i$ is the indicator random variable for event $A_i$. For indices $i,j$ we write $i\sim j$ if $i\neq j$ and $A_i, A_j$ are not independent. Set $$\Delta=\sum_{i\sim j}P[A_i\cap A_j],$$
    and $\mu=\textbf{E}X$. If $\Delta\ge \mu$, then 
    $$P[X=0]\le \exp\left(-\mu^2/2\Delta\right).$$}
    \begin{proof}
    This is Theorem 8.1.2 of~\cite{Chernoff from K and S}, adjusted to our settings.
    \end{proof}

    Using this, we can now show the (exponentially)-likely existence of a long cycle in the largest component:
    \begin{proof}[\textbf{\textit{Proof of Theorem 5}}]
    Let $p_1=\frac{1+\epsilon/2}{d}$, and let $\beta_1,\beta_2$ be some positive constants (possibly depending on $\epsilon$) to be determined later. Then, by Theorem 2 of~\cite{site-percolation}, with probability at least $1-\exp\left(-\beta_1\alpha^2n/d\right)$ there exists a path of length $\frac{\epsilon^2n}{20d}$ in $G[V_{p_1}]$, which is in a large component in $G[V_{p_1}]$. Consider the first $\frac{\epsilon^2n}{50d}$ vertices of the path. Recalling that, by assumption, $\alpha\le \epsilon^3$, we may apply Property 2 of Lemma 2.4, and conclude that their neighbourhood in $G$, disjoint from the path, is with probability at least $1-\exp\left(-\beta_2\alpha^2n/d\right)$ of size at least
    \begin{align*}
        (1-2\alpha)n\left(1-\exp(-\epsilon^2/50)\right)-\frac{\epsilon^2n}{20d}&\ge (1-2\alpha)n\left(\frac{\epsilon^2}{50}-\epsilon^4\right)-\frac{\epsilon^2n}{20d}\\
        &\ge\frac{\epsilon^2n}{400},
    \end{align*}
    where we used inequality $(7)$ and the fact that $d\ge 3$. We denote this neighbourhood by $N_1$. Similarly, considering the last $\frac{\epsilon^2n}{50d}$ vertices of the path, with probability at least $1-\exp\left(-\beta_2\alpha^2n/d\right)$ their neighbourhood in $G$, disjoint from the path, is of size at least $\frac{\epsilon^2n}{400}$. We denote this neighbourhood by $N_2$. By the expander mixing lemma (Lemma 2.2),
    \begin{align*}
        e_G(N_1,N_2)&\ge \frac{\epsilon^4nd}{400^2}-\frac{\lambda\epsilon^2n}{400}\\
                  &\ge c\epsilon^4nd,
    \end{align*}
    for some small enough $c>0$ a constant.
    
    Let $p_2=\frac{\epsilon}{2d-\epsilon-2}$, and consider $V_{p_2}$. Note that the distributions of $G[V_p]$ and $G[V_{p_1}\cup V_{p_2}]$ are identical. Now, if any of the edges between $N_1$ and $N_2$ has both of its endpoints in $V_{p_2}$, this edge will close a cycle of the required length together with the middle $\frac{\epsilon^2n}{100d}$ vertices of the path. We will show that with probability $1-\exp\left(-\beta\frac{\alpha^2n}{d}\right)$ one of these edges belongs to $G[V_{p_2}]$.
    
    For that, let $X$ be the random variable representing the number of edges between $N_1$ and $N_2$ that belong to $G[V_{p_2}]$. Observe that the expander mixing lemma could have counted each edge at most twice (and all of them twice, if $N_1=N_2$), and we thus have (with the same probability) at least $c\epsilon^4nd/2$ distinct edges between $N_1$ and $N_2$ in $G$. The probability one such an edge is in $G[V_{p_2}]$ is $p_2^2$, representing the choice of both of its ends. Hence:
    \begin{align*}
    \textbf{E}X&\ge \frac{c\epsilon^4nd}{2}p_2^2\\
    &\ge \frac{c\epsilon^6n}{8d}.
    \end{align*}
    Using the notation of Lemma 6.5, let $A_i$ be the event that the edge $e_i$, between $N_1$ and $N_2$, is in $G[V_{p_2}]$. Note that $A_i$ and $A_j$ are not independent only if the edges $e_i$ and $e_j$ intersect at some vertex $v$, and in that case $P[A_i\cap A_j]=p_2^3$. As for the possible number of pairs of intersecting edges, we have $n$ vertices where they may intersect, and we have at most ${d\choose 2}$ possible pairs of intersecting edges for each given vertex, and in total at most $nd^2$ pairs of intersecting edges in $G$. As such,
    $$\Delta\le nd^2p_2^3\le \frac{2\epsilon^3n}{d}.$$
    On the other hand, $\textbf{E}X\le\epsilon^4n/d$ and $\Delta\ge \epsilon^4n/d$, and as such by Lemma 6.5:
    $$P[X=0]\le\exp\left(-\mu^2/2\Delta\right)\le \exp\left(-\epsilon^{10}n/d\right).$$
    We thus obtain that $X>0$ with probability at least $1-\exp\left(-\beta\alpha^2n/d\right)$ for some positive constant $\beta(\epsilon)$. Since $X$ can only take integer values this implies $X\ge 1$ with the same probability, completing the proof.
    \end{proof}
    
    The proof of Theorem 6 will be similar to that of Theorem 4.1 in~\cite{high-girth}. 
    \begin{proof}[\textbf{\textit{Proof of Theorem 6}}]
    Let $p_1=\frac{1+\epsilon-\alpha}{d}$ and let $\rho=\frac{p_1}{p}$. Consider $V_{p_1}$. Observe that the random graph $G_1\sim G[V_{p_1}]$ can also be obtained by drawing a random graph $G\sim G[V_p]$ and then retaining each vertex in $G$ with $\rho$ independently. 
    
    Denote by $L_1$ the set of vertices of the largest component in $G[V_p]$, and by $L_1'$ the set of vertices of the largest component in $G[V_{p_1}]$. Applying Theorem 2 to $G[V_p]$ and $G[V_{p_1}]$, we have with probability at least $1-\exp\left(-\beta'\frac{\alpha^2n}{d}\right)$ (for some $\beta'=\beta'(\epsilon)$) that $$\big|L_1\big|\le \frac{xn}{d}+\frac{7 \alpha n}{d},$$ and $$\big|L_1'\big| \ge \frac{x_1 n}{d}-\frac{7\alpha n}{d},$$ where $x_1=x_1(\epsilon_1)$ is as defined by $(2)$, with $\epsilon_1=\epsilon-\alpha$. Recall that $x=2\epsilon-\frac{2\epsilon^2}{3}+O(\epsilon^3)$, and similarly $x_1=2(\epsilon-\alpha)-\frac{2(\epsilon-\alpha)^2}{3}+O\left((\epsilon-\alpha)^3\right)$. Thus,
    \begin{align*}
        x_1&=x-2\alpha+\frac{4\alpha\epsilon}{3}-\frac{2\alpha^2}{\epsilon}+O\left(\epsilon^2\alpha+\epsilon\alpha^2\right)\\
        &\ge x-2\alpha,
    \end{align*}
    given that $\alpha<\epsilon^2$ and for small enough $\epsilon$. Therefore, with probability at least $1-\exp\left(-\beta'\frac{\alpha^2n}{d}\right)$,
    $$\big|L_1'\big|>\frac{xn}{d}-\frac{9\alpha n}{d}.$$
    
    Let $\beta=\beta(\epsilon)$ be some positive constant to be determined later. Let $\mathcal{A}$ be the following event addressing $G[V_p]$:
    \begin{align*}
        \mathcal{A}=\Bigg\{\big|L_1\big|\le \frac{(x+7\alpha)n}{d} \text{ \& } \left(\exists S\subseteq L_1\text{, } \frac{16\alpha n}{d}\le |S|\le \frac{(x-10\alpha)n}{d}\text{, } \big|N_{G[V_p]}(S)\big|\le \frac{\beta\alpha^2}{\ln\left(\frac{1}{\alpha}\right)}\cdot\frac{n}{d}\right)\Bigg\}.
    \end{align*}
    Let $\mathcal{B}$ be the following event addressing $G[V_{p_1}]$:
    \begin{align*}
        \mathcal{B}=\Bigg\{\big|L_1'\big|\le \frac{(x-9\alpha)n}{d}\Bigg\}.
    \end{align*}
    
    Suppose $G[V_p]$ satisfies $\mathcal{A}$, and choose $S$ as in the definition of $\mathcal{A}$. If we remove all the vertices that neighbour $S$ in $G[V_p]$, we will separate $S$ from the rest of $L_1$. Since $\big|N_{G[V_p]}(S)\big|\le\frac{\beta\alpha^2}{\ln\left(\frac{1}{\alpha}\right)}\cdot\frac{n}{d}$, the probability to erase all these vertices when going from $G[V_p]$ to $G[V_{p_1}]$ is at least $(1-\rho)^{\frac{\beta\alpha^2}{\ln\left(\frac{1}{\alpha}\right)}\cdot\frac{n}{d}}$. But then, since  all other components in $G[V_p]$ are typically much smaller (and stay such in $G[V_{p_1}]$), we have that: $$\big|L_1'\big|\le \max\left(|S|, |L_1|-|S|\right)\le \frac{(x-9\alpha)n}{d}.$$
    Therefore,
    $$P[\mathcal{B}|\mathcal{A}]\ge (1-\rho)^{\frac{\beta\alpha^2}{\ln\left(\frac{1}{\alpha}\right)}\cdot\frac{n}{d}},$$
    and by Theorem 2, $$P[\mathcal{B}]\le\exp\left(-\beta'\frac{\alpha^2 n}{d}\right).$$
    We thus conclude that:
    \begin{align*}
    P[\mathcal{A}]&\le \frac{\exp\left(-\beta'\frac{\alpha^2 n}{d}\right)}{(1-\rho)^{\frac{\beta\alpha^2}{\ln\left(\frac{1}{\alpha}\right)}\cdot\frac{n}{d}}}\\
    &=\frac{\exp\left(-\beta'\frac{\alpha^2 n}{d}\right)}{\left(\frac{\alpha}{1+\epsilon}\right)^{\frac{\beta\alpha^2}{\ln\left(\frac{1}{\alpha}\right)}\cdot\frac{n}{d}}}\\
    &\le \exp\left(\left(2\beta-\beta'\right)\frac{\alpha ^2 n}{d}\right)= o(1),
    \end{align*}
    for $\beta<\frac{\beta'}{2}$ small enough. All that is left is to observe that the probability that an event violates the statement of this Theorem is at most $$P[\mathcal{A}]+\exp\left(-\beta'\frac{\alpha^2 n}{d}\right)+\left(\frac{n}{d}\right)^{-\beta''\ln\ln(n/d)}= o(1).$$
    \end{proof}

\paragraph{Acknowledgements.} The authors wish to thank Asaf Nachmias and Wojciech Samotij for their remarks on earlier versions of this paper, and to the anonymous referees for their careful reading and helpful remarks.

\end{document}